\documentclass[a4paper,12pt]{article} 

\usepackage{color}

\usepackage{amsfonts, amsmath, amsthm, amssymb}
\usepackage[T1]{fontenc}
\usepackage[cp1250]{inputenc}

\usepackage{graphicx}
\usepackage{amssymb}
\usepackage{amsmath}
\usepackage{mathptmx}
\usepackage{helvet}
\usepackage{courier}
\usepackage{txfonts}
\usepackage{tikz} 
   \usetikzlibrary{calc}
	\usetikzlibrary{arrows}
\usepackage{type1cm}

\usepackage{verbatim}

\usepackage{graphicx}
\usepackage{epsfig,amscd,amssymb,amsxtra,amsmath,amsthm}
\usepackage{type1cm}
\usepackage[T1]{fontenc}
\usepackage{graphics}
\usepackage[mathscr]{eucal}
\usepackage[all]{xy}
\usepackage{amsmath,amscd}

\newtheorem{theorem}{Theorem}[section]

\newtheorem{definition}[theorem]{Definition}
\newtheorem{lemma}[theorem]{Lemma}

\newtheorem{corollary}[theorem]{Corollary}
\newtheorem{problem}[theorem]{Problem}
\newtheorem{observation}[theorem]{Observation}

\newcommand{\diam} {\mathop{\rm diam}\nolimits}

\begin{document}

\def\joinrel{\mkern-3mu}
\newcommand{\varproj}{\displaystyle \lim_{\multimapinv\joinrel-\joinrel-}}

\title{Mapping and fixed point property theorems for inverse limits with set-valued bonding functions}
\author{Iztok Bani\v c,  Goran Erceg, and Judy Kennedy}
\date{}

\maketitle

\begin{abstract}
Among other results, the paper gives new mapping theorems and new fixed point property theorems for inverse limits of inverse sequences of compact metric spaces with upper semicontinuous set-valued bonding functions. We also revisit the results from two papers, Mioduszewski's ``Mappings of inverse limits'' \cite{mioduszewski} and Feuerbacher's ``Mappings of inverse limits revisited'' \cite{feuerbacher}.
\end{abstract}
\-
\\
\noindent
{\it Keywords:} continua, inverse limits, mapping theorem, fixed point property\\
\noindent
{\it 2020 Mathematics Subject Classification:} 54F17, 54H25, 37C25, 37B45, 54C60, 54F15, 37B45


\section{Introduction}\label{s1}
The motivation for this paper may be found in the following three areas of the inverse limit theory:
\begin{enumerate}
    \item Mapping theorems for inverse limits: our work is motivated by the Mioduszewski and Feuerbacher mapping theorems that are obtained in \cite{mioduszewski} and \cite{feuerbacher}.
    \item Fixed point property of inverse limits: our work is also motivated by the Feuerbacher fixed point property theorems from \cite{feuerbacher}  that may be used to construct examples of continua without the fixed point property. 
    \item Inverse limits with set-valued bonding functions: these are a generalization of inverse limits with continuous bonding functions and there is not much known about the fixed point property of such inverse limits. There is also not much known about existence of continuous surjections from one such inverse limit to another. However, some useful mapping theorems  for these generalized inverse limits are listed in \cite{ingi,ingiingi}.
\end{enumerate}

 Since any continuum $K$ may be expressed as an inverse limit $\varprojlim\{X_\ell,f_\ell\}_{\ell=1}^{\infty}$ of an inverse sequence of a certain family of simpler spaces (i.e., connected polyhedra), given a function $g:\varprojlim\{X_\ell,f_\ell\}_{\ell=1}^{\infty}\rightarrow  \varprojlim\{Y_\ell,g_\ell\}_{\ell=1}^{\infty}$, it is natural to consider the impact of the function  $g$ on the factor spaces $X_\ell$ and $Y_\ell$ of the inverse limits. This leads to a definition of the transformation of $g$ into a countable family $\{T_{n,m}(g)  \ | \  n,m \textup{ are positive integers}\}$ of upper semicontinuous set-valued functions $T_{n,m}(g):X_m\multimap Y_n$ on the factor spaces, and we exploit the relationship between these functions and the function $g$ to obtain new results about spaces with the fixed point property.  This idea has already been used by many authors to obtain interesting results about inverse limit spaces, for examples see \cite{loncar,mccord,rochowski,overover,rosen}.  We demonstrate the power of such a transformation of $g$ by applying it to obtain new results in the following categories.
 \begin{enumerate}
    \item Using the above-mentioned transformation of $g$, we introduce and prove new mapping theorems for inverse limits. They give a new insight into Mioduszewski's  and Feuerbacher's results from \cite{mioduszewski} and \cite{feuerbacher}. Among other things, we obtain necessary and sufficient conditions about existence of a continuous surjection from one inverse limit to another, see Corollary \ref{mappingTHM} and Theorem \ref{mappingTHM2}.
    \item Using transformations  $\{T_{n,m}(g)  \ | \  n,m \textup{ are positive integers}\}$, we obtain characterizations of inverse limits that do not have the fixed point property, for example, see Theorems \ref{Goran} and \ref{Goran1} and Corollary \ref{goran2}. Among other things, we revisit Feuerbacher's results from \cite{feuerbacher} to obtain another tool for constructing an inverse limit that does not have the fixed point property, see Theorem \ref{main:1} for details. Our method  allows us to work with upper semicontinuous set-valued functions between factor spaces in the construction of the almost commutative diagram that gives the desired inverse limit (see Figure \ref{diagram1})
\begin{figure}[h]
\centering
\noindent \begin{tikzpicture}[node distance=1.5cm, auto]
  \node (X1) {$X_{m_1}$};
  \node (X2) [right of=X1] {};
   \node (X3) [right of=X2] {$X_{m_2}$};
    \node (X4) [right of=X3] {};
      \node (X5) [right of=X4] {$X_{m_3}$};
    \node (X6) [right of=X5] {$\ldots$};
        \node (X7) [right of=X6] {$\varprojlim\{X_\ell,f_\ell\}_{\ell=1}^{\infty}$};
   \node (Z1) [below of=X1] {};
  \draw[<-] (X1) to node {$f_{m_1,m_2}$} (X3);
   \draw[<-] (X3) to node {$f_{m_2,m_3}$} (X5);
      \draw[<-] (X5) to node {$f_{m_3}$} (X6);
\node (Y1) [below of=Z1] {$Y_{n_1}$};
  \node (Y2) [right of=Y1] {};
    \node (Y3) [right of=Y2] {$Y_{n_2}$};
      \node (Y4) [right of=Y3] {};
       \node (Y5) [right of=Y4] {$Y_{n_3}$};
     \node (Y6) [right of=Y5] {$\ldots$};
      \node (Y7) [right of=Y6] {$\varprojlim\{Y_\ell,g_\ell\}_{\ell=1}^{\infty}$};
  \draw[<-] (Y1) to node {$g_{n_1,n_2}$} (Y3);
   \draw[<-] (Y3) to node {$g_{n_2,n_3}$} (Y5);
     \draw[<-] (Y5) to node {$g_{n_3}$} (Y6);
      \draw[o-] (Y1) to node {$H_1$} (X1);
            \draw[o-] (Y3) to node {$H_2$} (X3);
                  \draw[o-] (Y5) to node {$H_3$} (X5);
  \draw[<-] (Y7) to node {$g$} (X7);
\end{tikzpicture}
\caption{An almost commutative diagram with upper semicontinuous set-valued functions} \label{diagram1}
\end{figure}
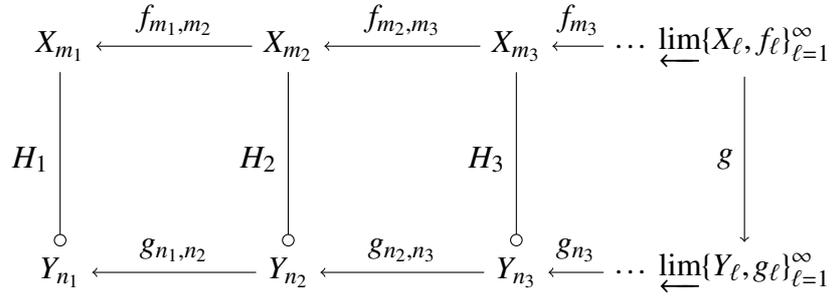
while the Feuerbacher method is restricted only to single-valued continuous functions (see Figure \ref{diagram2}).

 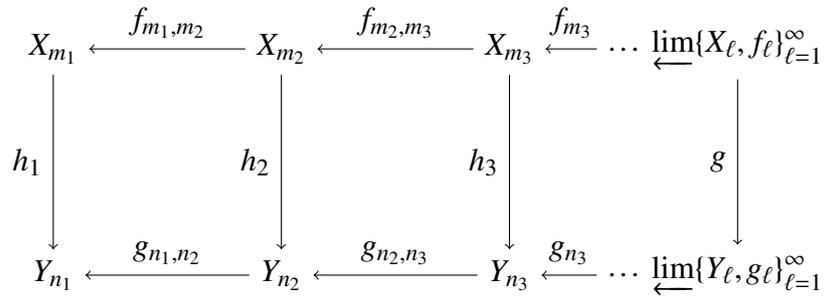
\begin{figure}[h]
\centering
\noindent \begin{tikzpicture}[node distance=1.5cm, auto]
  \node (X1) {$X_{m_1}$};
  \node (X2) [right of=X1] {};
   \node (X3) [right of=X2] {$X_{m_2}$};
    \node (X4) [right of=X3] {};
      \node (X5) [right of=X4] {$X_{m_3}$};
    \node (X6) [right of=X5] {$\ldots$};
        \node (X7) [right of=X6] {$\varprojlim\{X_\ell,f_\ell\}_{\ell=1}^{\infty}$};
   \node (Z1) [below of=X1] {};
  \draw[<-] (X1) to node {$f_{m_1,m_2}$} (X3);
   \draw[<-] (X3) to node {$f_{m_2,m_3}$} (X5);
      \draw[<-] (X5) to node {$f_{m_3}$} (X6);
\node (Y1) [below of=Z1] {$Y_{n_1}$};
  \node (Y2) [right of=Y1] {};
    \node (Y3) [right of=Y2] {$Y_{n_2}$};
      \node (Y4) [right of=Y3] {};
       \node (Y5) [right of=Y4] {$Y_{n_3}$};
     \node (Y6) [right of=Y5] {$\ldots$};
      \node (Y7) [right of=Y6] {$\varprojlim\{Y_\ell,g_\ell\}_{\ell=1}^{\infty}$};
  \draw[<-] (Y1) to node {$g_{n_1,n_2}$} (Y3);
   \draw[<-] (Y3) to node {$g_{n_2,n_3}$} (Y5);
     \draw[<-] (Y5) to node {$g_{n_3}$} (Y6);
      \draw[<-] (Y1) to node {$h_1$} (X1);
            \draw[<-] (Y3) to node {$h_2$} (X3);
                  \draw[<-] (Y5) to node {$h_3$} (X5);
                        \draw[<-] (Y7) to node {$g$} (X7);
\end{tikzpicture}
\caption{An almost commutative diagram with continuous single-valued functions} \label{diagram2}
\end{figure}

Also, our results are not restricted to inverse limits of polyhedra\footnote{The definition of a polyhedron may be found in \cite[p. 470--473]{poli}.}, they are obtained more generally, for inverse limits of arbitrary metric compacta. 
    \item We also apply the results about inverse limits with continuous single-valued functions to inverse limits of upper semicontinuous set-valued functions. By using these results we answer an open problem that is stated in \cite{ingi} by  W.~T.~Ingram  concerning fixed points and inverse limits with upper semicontinuous bonding functions.

\begin{problem}\cite[p.\ 81, Problem 6.53.\ and Problem 6.54.]{ingi}\label{problemcek}
If $(F_\ell)$ is a sequence of upper semicontinuous functions $F_\ell$ from $[0,1]$ into $2^{[0,1]}$ (Problem 6.53.) or into $C[0,1]$ (Problem 6.54.), what can be said about the fixed point property for the inverse limit $\varproj\{[0,1],F_\ell\}_{\ell=1}^{\infty}$?
\end{problem}
M.~M.~Marsh gave some very interesting results in  \cite{marsh1,marsh2} about the fixed point property of inverse limits of inverse sequences of closed intervals with upper semicontinuous set-valued bonding functions that give an answer to Problem \ref{problemcek}. In a way, our results are more general - they are given for inverse limits of any compact metric spaces, not only for the inverse limits of closed intervals. Among other things, our results give necessary and sufficient conditions for an inverse limit with upper semicontinuous set-valued bonding functions to have the fixed point property, see Theorem \ref{tralala}. We also give and prove new useful mapping theorems for inverse limits with upper semicontinuous set-valued bonding functions, i.e., we give necessary and sufficient conditions  about existence of a continuous surjection from one such inverse limit to another, see Theorem \ref{hopsasa}.
\end{enumerate}

We proceed as follows. In Section \ref{s2}, basic definitions and well-known results that are needed in the paper are given, in Section \ref{s3}, we introduce the concept of $T_{n,m}$-transformations of functions $g:\varprojlim\{X_\ell,f_\ell\}_{\ell=1}^{\infty}\rightarrow  \varprojlim\{Y_\ell,g_\ell\}_{\ell=1}^{\infty}$ and give new results about the fixed point property of inverse limits. In Section \ref{s4}, we apply the $T_{n,m}$-transformations to generalize Mioduszewski's  and Feuerbacher's results from \cite{mioduszewski} and \cite{feuerbacher}. Finally, in Section \ref{s5}, we apply the obtained results to inverse limits with upper semicontinuous set-valued bonding functions.

\section{Definitions and notation}\label{s2}

We use $\mathbb N$ to denote the set of positive integers.
All spaces in this paper are metric.

\begin{definition}
A {\em fixed point\/} of a function $f:X\rightarrow X$ is a point $p$ in $X$ such that $f(p)=p$. A topological space $X$ is said to have the {\em fixed point property\/} if every continuous function from $X$ into $X$ has a fixed point.
\end{definition}

\begin{definition}
A {\em continuum\/}  is a non-empty connected compact metric space. A continuum is {\em degenerate\/} if it consists of only a single point. Otherwise it is {\em non-degenerate\/}.  A {\em subcontinuum\/} is a subspace of a continuum which itself is also a continuum.
\end{definition}

\begin{definition}
Let $\ell$ be a positive integer,  $X$ a metric space, and $f:X\rightarrow X$ a function. We use $f^\ell$ to denote the {\em composition\/}: 
$$
f^\ell=\underbrace{f\circ f\circ\cdots\circ f}_\ell. 
$$
\end{definition}

\begin{definition}
Let $X$ be a metric space. The {\em identity function\/} on $X$ will be denoted by $1_{X}$.
\end{definition}

\begin{definition}
For each positive integer $\ell$, let $(X_\ell,d_\ell)$ be a compact metric space, where $d_{\ell}(x,y)<1$ for all $x,y\in X_{\ell}$. We use $D$ to denote the {\em product metric\/} on $\prod_{\ell=1}^{\infty}X_\ell$, defined by 
$$
D(\mathbf{x}, \mathbf{y})=\sup\left\{\frac{d_\ell(x_\ell,y_\ell)}{2^\ell} \ | \ \ell \textup{ is a positive integer }\right\}
$$ 
for all $\mathbf{x}=(x_1,x_2,x_3,\ldots), \mathbf{y}=(y_1,y_2,y_3,\ldots)\in \prod_{\ell=1}^{\infty}X_\ell$.  
\end{definition}
\begin{definition}
If $(X,d)$ is a compact metric space, then $2^X$ denotes the set of all {\em nonempty closed subsets\/} of $X$, and $C(X)$ the set of all {\em connected elements\/} of $2^X$.  
\end{definition}
\begin{definition}
Let $X$ and $Y$ be compact metric spaces. 
A function $F: X\rightarrow 2^Y$ is called a {\em set-valued function\/} from $X$ to $Y$. We denote set-valued functions $F: X\rightarrow 2^Y$ by $F: X\multimap Y$.
\end{definition}
\begin{definition}
A function  $F : X\multimap Y$
is an {\em upper semicontinuous\/} set-valued function if for each
open set  $V\subseteq Y$ the set $\{x\in X \ | \ F(x) \subseteq V\}$ is an
open set in $X$.
\end{definition}
\begin{definition}
The {\em graph\/} $\Gamma(F)$ of a set-valued function $F:X\multimap Y$ is the set of
all points  $(x,y)\in X\times Y$ such that $y \in F(x)$. 
\end{definition}
There is a simple characterization of upper semicontinuous set-valued functions
(\cite[Theorem 1.2, p.\ 3]{ingi}):

\begin{theorem}
\label{th:grafi}  Let $X$ and $Y$ be compact metric spaces and $F:X\multimap Y$ a set-valued function. Then $F$ is upper semicontinuous if and only if its
graph $\Gamma(F)$ is closed in  $X\times Y$. 
\end{theorem}

If $F:X\multimap Y$ is a set-valued function, where for each $x \in X$, the image $F(x)$ is a singleton in $Y$, then we can interpret it as a 
single-valued function, identifying it with the function $f:X\rightarrow Y$, where $F(x)=\{f(x)\}$ for any $x\in X$. Conversely, any single-valued function $f:X\rightarrow Y$ can be identified with the set-valued function $F:X\multimap Y$, defined by $F(x)=\{f(x)\}$. 
\begin{definition}
 Let $F:X\multimap Y$ be a set-valued function. Then we define $F^{-1}:F(X)=\bigcup_{x\in X}F(x)\multimap X$ by $F^{-1}(y)=\{x\in X \ | \ y\in F(x)\}$ for any $y\in Y$.
 \end{definition}
\begin{definition}
We say that the graph of a set-valued function $F:X\multimap Y$ is 
{\em surjective\/} if for each $y\in Y$, $|F^{-1}(y)|\geq 1$, i.e.\ if $F(X)=Y$.
\end{definition}

In this paper we mostly deal with {\em inverse sequences\/} $\{X_\ell,F_\ell\}_{\ell=1}^{\infty}$,  
where $X_\ell$ are compact metric spaces and $F_\ell:X_{\ell+1}\multimap {X_\ell}$ upper semicontinuous functions. The spaces $X_\ell$ are called {\em factor spaces\/} and the functions $F_\ell$ are called the {\em bonding functions\/}.
\begin{definition}
The {\em inverse limit\/} of an inverse sequence $\{X_\ell,F_\ell\}_{\ell=1}^\infty$ is defined to be the  subspace of the product space $\prod_{\ell=1}^\infty X_\ell$ of all $\mathbf x=(x_1,x_2,x_3,\ldots)\in \prod_{\ell=1}^\infty X_\ell$,  such that $x_\ell \in F_\ell(x_{\ell+1})$ for each $\ell$. The inverse limit is denoted by  $\varproj\{X_\ell,F_\ell\}_{\ell=1} ^\infty$.
\end{definition}
These inverse limits are a generalization (introduced by T.~W.~Ingram and W.~S.~Mahavier) of inverse limits of inverse sequences $\{X_\ell,f_\ell\}_{\ell=1}^{\infty}$,  
where $X_\ell$ are compact metric spaces and $f_\ell:X_{\ell+1}\rightarrow {X_\ell}$ continuous functions. Such inverse limits are usually denoted by $\varprojlim\{X_\ell,f_\ell\}_{i=1} ^\infty$. Obviously, for any inverse sequence $\{X_\ell,f_\ell\}_{\ell=1}^{\infty}$ of compact metric spaces and continuous bonding functions, 
$$
\varprojlim\{X_\ell,f_\ell\}_{\ell=1} ^\infty=\varproj\{X_\ell,F_\ell\}_{\ell=1} ^\infty,
$$
where $F_{\ell}:X_{\ell+1}\multimap X_\ell$ is defined by $F_{\ell}(x)=\{f_{\ell}(x)\}$ for every $x\in X_{\ell}$.

\section{$T_{n,m}$-transformations}\label{s3}
In the present section, we define our main tool that is used in the paper, the $T_{n,m}$-transformations of continuous functions from one inverse limit to another, and prove some basic results about such transformations. Before we do that, we define special kinds of projections in the following definition. 
\begin{definition}
Let $\{X_\ell,f_\ell\}_{\ell=1}^{\infty}$ and $\{Y_\ell,g_\ell\}_{\ell=1}^{\infty}$ be inverse sequences of compact metric spaces and surjective continuous bonding functions. We use $p_n$ to denote the {\em $n$-th standard projection\/} $p_n:\varprojlim\{X_\ell,f_\ell\}_{\ell=1}^{\infty}\rightarrow X_n$, and $q_n$ to denote the {\em $n$-th standard projection\/} $q_n:\varprojlim\{Y_\ell,g_\ell\}_{\ell=1}^{\infty}\rightarrow Y_n$.
\end{definition}

\begin{definition}
Let $\{X_\ell,f_\ell\}_{\ell=1}^{\infty}$ and $\{Y_\ell,g_\ell\}_{\ell=1}^{\infty}$ be inverse sequences of compact metric spaces and surjective continuous bonding functions. For any function $g: \varprojlim\{X_\ell,f_\ell\}_{\ell=1}^{\infty} \rightarrow \varprojlim\{Y_\ell,g_\ell\}_{\ell=1}^{\infty}$ and for any pair $(n,m)$ of positive integers we use $T_{n,m}(g)$ to denote the {\em $T_{n,m}$-transformation\/} $T_{n,m}(g):X_m\multimap Y_n$ of the function $g$, which is defined by
$$
T_{n,m}(g)(x)=q_n(g(p_m^{-1}(x)))
$$
for any $x\in X_m$ (se Figure \ref{diagram3}).

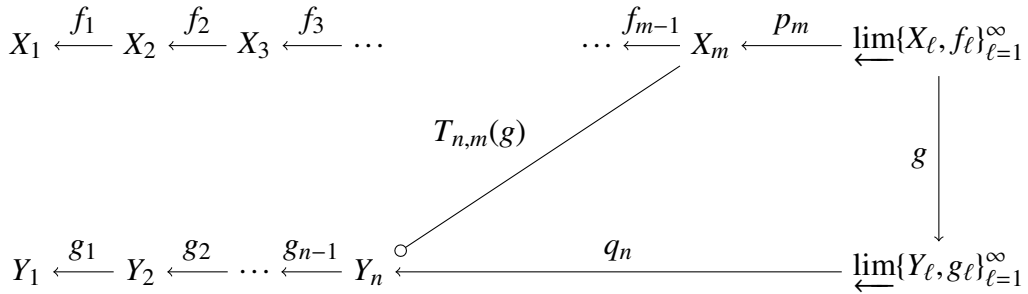
\begin{figure}[h]
\centering
\noindent\begin{tikzpicture}[node distance=1.5cm, auto]
  \node (X1) {$X_1$};
  \node (X2) [right of=X1] {$X_{2}$};
    \node (X3) [right of=X2] {$X_{3}$};
    \node (X4) [right of=X3] {$\ldots$};
      \node (X5) [right of=X4] {};
    \node (X6) [right of=X5] {$\ldots$};
        \node (X7) [right of=X6] {$X_{m}$};
    \node (X8) [right of=X7] {};
    \node (X9) [right of=X8] {$\varprojlim\{X_\ell,f_\ell\}_{\ell=1}^{\infty}$};
   \node (Z1) [below of=X1] {};
  \node (Z2) [right of=Z1] {};
    \node (Z3) [right of=Z2] {};
    \node (Z4) [right of=Z3] {};
      \node (Z5) [right of=Z4] {};
    \node (Z6) [right of=Z5] {};
        \node (Z7) [right of=Z6] {};
    \node (Z8) [right of=Z7] {};
    \node (Z9) [right of=Z8] {};
  \draw[<-] (X1) to node {$f_1$} (X2);
   \draw[<-] (X2) to node {$f_2$} (X3);
     \draw[<-] (X3) to node {$f_3$} (X4);
         \draw[<-] (X6) to node {$f_{m-1}$} (X7);
          \draw[<-] (X7) to node {$p_{m}$} (X9);
\node (Y1) [below of=Z1] {$Y_1$};
  \node (Y2) [right of=Y1] {$Y_{2}$};
    \node (Y3) [right of=Y2] {$\ldots$};
     \node (Y4) [right of=Y3] {$Y_n$};
      \node (Y5) [right of=Y4] {};
      \node (Y6) [below of=Z9] {$\varprojlim\{Y_\ell,g_\ell\}_{\ell=1}^{\infty}$};
  \draw[<-] (Y1) to node {$g_1$} (Y2);
   \draw[<-] (Y2) to node {$g_2$} (Y3);
     \draw[<-] (Y3) to node {$g_{n-1}$} (Y4);
      \draw[<-] (Y4) to node {$q_n$} (Y6);
      \draw[o-] (Y4) to node {$T_{n,m}(g)$} (X7);
       \draw[<-] (Y6) to node {$g$} (X9);
\end{tikzpicture}
\caption{The set-valued function $T_{n,m}(g)$} \label{diagram3}
\end{figure}
\end{definition}
\begin{observation}
Let $\{X_\ell,f_\ell\}_{\ell=1}^{\infty}$ and $\{Y_\ell,g_\ell\}_{\ell=1}^{\infty}$ be inverse sequences of compact metric spaces and surjective continuous bonding functions. For any function $g: \varprojlim\{X_\ell,f_\ell\}_{\ell=1}^{\infty} \rightarrow \varprojlim\{Y_\ell,g_\ell\}_{\ell=1}^{\infty}$ and for any pair $(n,m)$ of positive integers, the function $T_{n,m}(g):X_m\multimap Y_n$ is an upper semicontinuous set-valued function whose graph is surjective and connected. 
\end{observation}
\begin{definition}
We use the following notation
$$
\mathcal W=\{(n,m)\in \mathbb N\times \mathbb N \ | \ n\leq m\}.
$$
\end{definition}
We prove the following lemma before we give and prove our characterization theorems for inverse limits that do not have the fixed point property. 
\begin{lemma}\label{Iztok}
Let $\{X_\ell,f_\ell\}_{\ell=1}^{\infty}$ be an inverse sequence of compact metric spaces and surjective continuous bonding functions, and let $g: \varprojlim\{X_\ell,f_\ell\}_{\ell=1}^{\infty} \rightarrow \varprojlim\{X_\ell,f_\ell\}_{\ell=1}^{\infty}$ be a continuous function.  Also, let $(n_0,m_0)\in \mathcal W$ be such that $$
f_{n_0,m_0}(x)\not \in T_{n_0,m_0}(g)(x)
$$
for each $x\in X_{m_0}$. Then for all $(n,m)\in \mathcal W$, the following holds
$$
n\geq n_0, m\geq m_0 \Longrightarrow f_{n,m}(x)\not \in T_{n,m}(g)(x) \textup{ for each } x\in X_m.
$$
\end{lemma}
\begin{proof}
Let $(n,m)\in \mathcal W$ such that $n\geq n_0$ and $m\geq m_0$. 
Suppose that there is a point $x_0\in X_{m}$ such that $f_{n,m}(x_0) \in T_{n,m}(g)(x_0)$. Fix such an element $x_0$ and let $y_0=f_{m_0,m}(x_0)$. We show that $f_{n_0,m_0}(y_0)\in T_{n_0,m_0}(g)(y_0)$.

Let $\mathbf z=(z_1,z_2,z_3,\ldots)\in p_{m}^{-1}(x_0)$ such that $p_{n}(g(\mathbf z))=f_{n,m}(x_0)$.  Such a point $\mathbf z$ does exist since
$$
f_{n,m}(x_0) \in T_{n,m}(g)(x_0) ~~~ \textup{ and } ~~~  T_{n,m}(g)(x_0)=p_{n}(g(p_{m}^{-1}(x_0))).
$$
Therefore,
$$
p_{n_0}(g(\mathbf z))=f_{n_0,n}(p_{n}(g(\mathbf z)))=f_{n_0,n}(f_{n,m}(x_0))
$$
$$
=f_{n_0,m}(x_0) =f_{n_0,m_0}(f_{m_0,m}(x_0))=f_{n_0,m_0}(y_0).
$$
follows. We claim that $\mathbf z\in p_{m_0}^{-1}(y_0)$. To prove the claim, recall that  $\mathbf z\in p_{m}^{-1}(x_0)$. Then $p_{m}(\mathbf z)=x_0$ and, therefore, $f_{m_0,m}(p_{m}(\mathbf z))=f_{m_0,m}(x_0)$. It follows that 
$$
p_{m_0}(\mathbf z)=f_{m_0,m}(p_{m}(\mathbf z))=f_{m_0,m}(x_0)=y_0
$$  
and we have proved the claim.
Therefore,
$$
f_{n_0,m_0}(y_0)=p_{n_0}(g(\mathbf z))\in p_{n_0}(g(p_{m_0}^{-1}(y_0))) = T_{n_0,m_0}(g)(y_0)
$$
follows.
\end{proof}

\begin{theorem}\label{Goran}
		Let $\{X_\ell,f_\ell\}_{\ell=1}^{\infty}$ be an inverse sequence of compact metric spaces and surjective continuous bonding functions,  and let $g: \varprojlim\{X_\ell,f_\ell\}_{\ell=1}^{\infty} \rightarrow \varprojlim\{X_\ell,f_\ell\}_{\ell=1}^{\infty}$ be a continuous function. The following statements are equivalent. 
			\begin{enumerate}
			\item \label{ena} For each $\mathbf x\in \varprojlim\{X_\ell,f_\ell\}_{\ell=1}^{\infty}$, $g(\mathbf x)\neq \mathbf x$.
			\item \label{dva} There is $(n,m)\in \mathcal W$ such that $f_{n,m}(x)\not \in T_{n,m}(g)(x)$ for each $x\in X_m$.
			\item \label{tri} There is a positive integer $N$ such that for each $(n,m)\in \mathcal W$,  it holds that 
			$$
			m\geq n\geq N \Longrightarrow f_{n,m}(x)\not \in T_{n,m}(g)(x) \textup{ for each } x\in X_m.
			$$
			\end{enumerate}
			\end{theorem}
\begin{proof}
First we prove the implication from \ref{ena} to \ref{dva}. Suppose that for each $(n,m)\in \mathcal W$, there is a point $x\in X_m$ such that $f_{n,m}(x)\in T_{n,m}(g)(x)$. So, for each positive integer $m$, let $x_m\in X_m$ be such a point that 
$$
f_{m-1,m}(x_m)\in T_{m-1,m}(g)(x_m).
$$
Let $m$ be any positive integer such that $m\geq 2$. 
Since $T_{m-1,m}(g)(x_m)=p_{m-1}(g(p_m^{-1}(x_m)))$, there is $\mathbf{z}^m=(z_1^m,z_2^m,z_3^m,\ldots )\in p_m^{-1}(x_m)$ (note that all the bonding functions are surjective) such that
$$
f_{m-1,m}(x_m)=p_{m-1}(g(\mathbf{z}^m)).
$$
Fix such a point $\mathbf{z}^m$.  It follows from $p_m(\mathbf{z}^m)=x_m$ and $p_{m-1}=f_{m-1}\circ p_m$ that 
$$
p_{m-1}(\mathbf{z}^m)=f_{m-1}(p_m(\mathbf{z}^m))=f_{m-1,m}(x_m)=p_{m-1}(g(\mathbf{z}^m)),
$$
and from this it follows that
$$
p_{m-2}(\mathbf{z}^m)=f_{m-2}(p_{m-1}(\mathbf{z}^m))=f_{m-2}(p_{m-1}(g(\mathbf{z}^m)))
=p_{m-2}(g(\mathbf{z}^m)).
$$
Inductively, we get that for any $k<m$, 
$$
p_k(\mathbf{z}^m)=p_k(g(\mathbf{z}^m)).
$$
Let $(\mathbf{z}^{i_m})$ be a convergent subsequence of the sequence $(\mathbf{z}^m)$ and let $\mathbf{z}^0=(z_1,z_2,z_3,\ldots )\in \varprojlim\{X_\ell,f_\ell\}_{\ell=1}^{\infty}$ such that $\displaystyle \lim_{m\to \infty}\mathbf{z}^{i_m}=\mathbf{z}^{0}$. Let $n$ be any positive integer. Then 
$$
p_n(g(\mathbf{z}^{0}))=p_n(g(\lim_{m\to \infty}\mathbf{z}^{i_m}))=\lim_{m\to \infty}p_n(g(\mathbf{z}^{i_m}))
$$
$$
=\lim_{m\to \infty}p_n(\mathbf{z}^{i_m})=p_n(\lim_{m\to \infty}\mathbf{z}^{i_m}))=p_n(\mathbf{z}^0).
$$
Therefore, $g(\mathbf{z}^0)=\mathbf{z}^0$, which is a contradiction.

Next we prove the implication from \ref{dva} to \ref{ena}. Let $(n,m)\in \mathcal W$ such that $f_{n,m}(x)\not \in T_{n,m}(g)(x)$ for each $x\in X_m$.

Suppose that there is $\mathbf x\in \varprojlim\{X_\ell,f_\ell\}_{\ell=1}^{\infty}$ such that $g(\mathbf x)= \mathbf x$.  Fix such an $\mathbf x=(x_1,x_2,x_3,\ldots)$.  Then 
$$
\mathbf x\in p_m^{-1}(x_m)
$$
and therefore, 
$$
g(\mathbf x)\in g(p_m^{-1}(x_m)).
$$
Since $g(\mathbf x)= \mathbf x$, it follows that
$$
\mathbf x\in g(p_m^{-1}(x_m)).
$$
Therefore,
$$
p_n(\mathbf x)\in p_n(g(p_m^{-1}(x_m))). 
$$
Since $f_{n,m}(x_m)=x_n=p_n(\mathbf x)$ and $p_n(g(p_m^{-1}(x_m)))=T_{n,m}(g)(x_m)$, it follows that 
$$
f_{n,m}(x_m)\in T_{n,m}(g)(x_m)
$$
which is a contradiction. 

The implication from \ref{tri} to \ref{dva} is obvious.  To complete the proof, we prove the implication from \ref{dva} to \ref{tri}. 
		
Let $(n_0,m_0)\in \mathcal W$ be such that $f_{n_0,m_0}(x)\not \in T_{n_0,m_0}(g)(x)$ for each $x\in X_{m_0}$ and let $N=m_0$. Let $(n,m)\in \mathcal W$ be any point such that $n\geq N$.  Since $n_0\leq n$ and $m_0\leq m$, it follows from Lemma \ref{Iztok} that $f_{n,m}(x)\not \in T_{n,m}(g)(x)$  for each $x\in X_m$.
\end{proof}
In the following theorem, we give two characterizations of  continua $X$  that do not have the fixed point property.

\begin{theorem}\label{Goran1}
Let $X$ be any continuum.  The following statements are equivalent.
\begin{enumerate}
\item \label{one}$X$ does not have the fixed point property.
\item \label{two} There is an inverse sequence $\{P_\ell,f_\ell\}_{\ell=1}^{\infty} $ of connected polyhedra $P_\ell$ and continuous surjections $f_\ell:P_{\ell+1}\rightarrow P_\ell$ such that $\varprojlim\{P_\ell,f_\ell\}_{\ell=1}^{\infty} $ is homeomorphic to $X$,  and there is a continuous function  
$$
g:\varprojlim\{P_\ell,f_\ell\}_{\ell=1}^{\infty}\rightarrow \varprojlim\{P_\ell,f_\ell\}_{\ell=1}^{\infty}$$
 such that there exists $(n,m)\in \mathcal W$ such that
			$$
		 f_{n,m}(x)\not \in T_{n,m}(g)(x) 
			$$
			holds for each $ x\in P_m$.
\item \label{three} There is an inverse sequence $\{P_\ell,f_\ell\}_{\ell=1}^{\infty} $ of connected polyhedra $P_\ell$ and continuous surjections $f_\ell:P_{\ell+1}\rightarrow P_\ell$ such that $\varprojlim\{P_\ell,f_\ell\}_{\ell=1}^{\infty} $ is homeomorphic to $X$,  and there is a continuous function  
$$
g:\varprojlim\{P_\ell,f_\ell\}_{\ell=1}^{\infty}\rightarrow \varprojlim\{P_\ell,f_\ell\}_{\ell=1}^{\infty}$$
 such that for each $(n,m)\in \mathcal W$,
			$$
		 f_{n,m}(x)\not \in T_{n,m}(g)(x)
			$$
			holds  for each $ x\in P_m$.
\end{enumerate}
\end{theorem}
\begin{proof}
Note that by \cite[Theorem 7, page 61]{mardesic} or \cite[Theorem 2.13, page 24]{nadler}, any continuum is homeomorphic to an inverse limit of inverse sequence of connected polyhedra with surjective continuous  bonding functions.  Suppose that $X$ does not have the fixed point property.  Let $\{P_\ell,f_\ell\}_{\ell=1}^{\infty} $ be an inverse sequence of connected polyhedra $P_\ell$ and continuous surjections $f_\ell:P_{\ell+1}\rightarrow P_\ell$ such that $\varprojlim\{P_\ell,f_\ell\}_{\ell=1}^{\infty} $ is homeomorphic to $X$.  Then $\varprojlim\{P_\ell,f_\ell\}_{\ell=1}^{\infty} $ does not have the fixed point property.  Let $g:\varprojlim\{P_\ell,f_\ell\}_{\ell=1}^{\infty} \rightarrow \varprojlim\{P_\ell,f_\ell\}_{\ell=1}^{\infty} $ be a continuous function that does not have fixed points.   By Theorem \ref{Goran},   there is $(n,m)\in \mathcal W$ such that
$$
		 f_{n,m}(x)\not \in T_{n,m}(g)(x) \textup{ for each } x\in P_m
			$$
	and we proved that \ref{two} follows from \ref{one}. To see that \ref{three} follows from \ref{two}, let  $\{P_\ell',f_\ell'\}_{\ell=1}^{\infty} $ be an inverse sequence of connected polyhedra $P_\ell'$ and continuous surjections $f_\ell':P_{\ell+1}'\rightarrow P_\ell'$ such that $\varprojlim\{P_\ell',f_\ell'\}_{\ell=1}^{\infty} $ is homeomorphic to $X$,  and let 
$$
g':\varprojlim\{P_\ell',f_\ell'\}_{\ell=1}^{\infty}\rightarrow \varprojlim\{P_\ell',f_\ell'\}_{\ell=1}^{\infty}$$
be a continuous function such that for some $(n,m)\in \mathcal W$,
			$$
		 f_{n,m}'(x)\not \in T_{n,m}(g')(x) \textup{ for each } x\in P_m'.
			$$

Let $N$ be a positive integer such that for each $(n,m)\in \mathcal W$,  it holds that 
			$$
			m\geq n\geq N \Longrightarrow f_{n,m}'(x)\not \in T_{n,m}(g')(x) \textup{ for each } x\in P_m'.
			$$ Such a positive integer $N$ does exist by Theorem \ref{Goran}.  
			For each positive integer $\ell$, let
			$$
			P_\ell=P_{N-1+\ell}'
			$$
			and
			$$
			f_\ell=f_{N-1+\ell}'.
			$$
			Also, let
			$$
			g=\sigma_{N}\circ g'\circ \sigma_N^{-1},
			$$
			where $\sigma_N: \varprojlim\{P_\ell',f_\ell'\}_{\ell=1}^{\infty} \rightarrow \varprojlim\{P_\ell',f_\ell'\}_{\ell=N}^{\infty} $ is the shift homeomorphism defined by
			$$
			\sigma_N(x_1,x_2,x_3,\ldots )=(x_N,x_{N+1},x_{N+2},\ldots ).
			$$
For each positive integer $m$, let  $p_m:\varprojlim\{P_\ell,f_\ell\}_{\ell=N}^{\infty}\rightarrow P_m$ be the $m$-th standard projection, defined by $p_m(\mathbf x)=x_m$ for each $\mathbf x=(x_1,x_2,x_3,\ldots)\in \varprojlim\{P_\ell,f_\ell\}_{\ell=N}^{\infty}$,  and  let $p_m':\varprojlim\{P_\ell',f_\ell'\}_{\ell=N}^{\infty}\rightarrow P_m'$ be the $m$-th standard projection, defined by $p_m'(\mathbf x)=x_m$ for each $\mathbf x=(x_1,x_2,x_3,\ldots)\in \varprojlim\{P_\ell',f_\ell'\}_{\ell=N}^{\infty}$.  Note that for any positive integer $\ell$,
$$
p_\ell\circ \sigma_N=p_{N-1+\ell}'.
$$
Let $(n,m)\in \mathcal W$. We show that
$ f_{n,m}(x)\not \in T_{n,m}(g)(x)$ for each $x\in P_m$.  Let $x\in P_m$. Then 
$$
f_{n,m}(x)=f_{N-1+n,N-1+m}'(x)
$$
and 
$$
T_{n,m}(g)(x)=p_n(g(p_m^{-1}(x)))=p_n(\sigma_N(g'(\sigma_N^{-1}(p_m^{-1}(x)))))
$$
$$
=p_{N-1+n}'(g'(p_{N-1+m}'^{-1}(x)))=T_{N-1+n,N-1+m}(g')(x).
$$
Since $N-1+n\geq N$, it follows that $f_{N-1+n,N-1+m}'(x)\not \in T_{N-1+n,N-1+m}(g')(x)$. Therefore, 
$$
f_{n,m}(x)\not \in T_{n,m}(g)(x).
$$
Finally we prove the implication from \ref{three}  to  \ref{one}.  Let  $\{P_\ell,f_\ell\}_{\ell=1}^{\infty} $ be an inverse sequence of connected polyhedra $P_\ell$ and continuous surjections $f_\ell:P_{\ell+1}\rightarrow P_\ell$ such that $\varprojlim\{P_\ell,f_\ell\}_{\ell=1}^{\infty} $ is homeomorphic to $X$,  and let 
$$
g:\varprojlim\{P_\ell,f_\ell\}_{\ell=1}^{\infty}\rightarrow \varprojlim\{P_\ell,f_\ell\}_{\ell=1}^{\infty}$$
be a continuous function such that for each $(n,m)\in \mathcal W$,  it holds that 
			$$
		 f_{n,m}(x)\not \in T_{n,m}(g)(x) \textup{ for each } x\in P_m.
			$$
	By Theorem \ref{Goran},   the function $g$ does not have any fixed points. Therefore, the inverse limit $\varprojlim\{P_\ell,f_\ell\}_{\ell=1}^{\infty} $ does not have the fixed point property.  Since $X$ is homeomorphic to $\varprojlim\{P_\ell,f_\ell\}_{\ell=1}^{\infty} $, also $X$ does not have the fixed point property. 
\end{proof}

In the following corollary, infinitely many characterizations of continua that do not have the fixed point property are obtained. 
\begin{corollary}\label{goran2}
Let $X$ be any continuum.  The following statements are equivalent for any pair $(n_0,m_0)\in \mathcal W$.
\begin{enumerate}
\item \label{one1} $X$ does not have the fixed point property.
\item \label{two2} There is an inverse sequence $\{P_\ell,f_\ell\}_{\ell=1}^{\infty} $ of connected polyhedra $P_\ell$ and continuous surjections $f_\ell:P_{\ell+1}\rightarrow P_\ell$ such that $\varprojlim\{P_\ell,f_\ell\}_{\ell=1}^{\infty} $ is homeomorphic to $X$,  and there is a continuous function  
$$
g:\varprojlim\{P_\ell,f_\ell\}_{\ell=1}^{\infty}\rightarrow \varprojlim\{P_\ell,f_\ell\}_{\ell=1}^{\infty}$$
 such that 
			$$
		 f_{n_0,m_0}(x)\not \in T_{n_0,m_0}(g)(x) \textup{ for each } x\in P_{m_0}.
			$$
			\end{enumerate}
\end{corollary}
\begin{proof}
The corollary follows directly from Theorem \ref{Goran1}.
\end{proof}
\section{Mioduszewski's and Feuerbacher's results revisited using $T_{n,m}$-transformations}\label{s4}
In this section, we revisit the Mioduszewski and Feuerbacher results from \cite{feuerbacher,mioduszewski} about mappings of inverse limits. Note that Mioduszewski and Feuerbacher's results are obtained for inverse sequences of polyhedra and then, it is discussed that the results can also be applied to inverse limits of inverse sequences of any metric compacta. For these particular inverse limits, they construct (also using tools from algebraic topology) continuous functions $h_{m,n}$ that have similar properties as our upper semicontinuous functions $T_{n,m}(g)$.  
\begin{theorem}\label{Mio:USC}
Let $\{X_\ell,f_\ell\}_{\ell=1}^{\infty}$ and $\{Y_\ell,g_\ell\}_{\ell=1}^{\infty}$ be inverse sequences of compact metric spaces and surjective continuous bonding functions, let $g:\varprojlim\{X_\ell,f_\ell\}_{\ell=1}^{\infty}\rightarrow \varprojlim\{Y_\ell,g_\ell\}_{\ell=1}^{\infty}$ be a continuous function, and let $\varepsilon>0$ and $n$ a positive integer.  Then there is a positive integer $m_0$ such that for each positive integer $m\geq m_0$,   
\begin{enumerate}
\item $\diam (T_{n,m}(g)(x))<\varepsilon $ for any $x\in X_m$,
\item $q_n(g(\mathbf x))\in  T_{n,m}(g)(p_m(\mathbf x))$ for each $\mathbf x\in\varprojlim\{X_\ell,f_\ell\}_{\ell=1}^{\infty}$.
\end{enumerate}
[See Figure \ref{diagram4}.]
\begin{figure}[h]
\centering
\noindent\begin{tikzpicture}[node distance=1.5cm, auto]
  \node (X1) {$X_1$};
  \node (X2) [right of=X1] {$X_{2}$};
    \node (X3) [right of=X2] {$X_{3}$};
    \node (X4) [right of=X3] {$\ldots$};
      \node (X5) [right of=X4] {$X_{m_0}$};
    \node (X6) [right of=X5] {$\ldots$};
        \node (X7) [right of=X6] {$X_{m}$};
    \node (X8) [right of=X7] {};
    \node (X9) [right of=X8] {$\varprojlim\{X_\ell,f_\ell\}_{\ell=1}^{\infty}$};
   \node (Z1) [below of=X1] {};
  \node (Z2) [right of=Z1] {};
    \node (Z3) [right of=Z2] {};
    \node (Z4) [right of=Z3] {};
      \node (Z5) [right of=Z4] {};
    \node (Z6) [right of=Z5] {};
        \node (Z7) [right of=Z6] {};
    \node (Z8) [right of=Z7] {};
    \node (Z9) [right of=Z8] {};
  \draw[<-] (X1) to node {$f_1$} (X2);
   \draw[<-] (X2) to node {$f_2$} (X3);
     \draw[<-] (X3) to node {$f_3$} (X4);
   \draw[<-] (X4) to node {$f_{m_0-1}$} (X5);
      \draw[<-] (X5) to node {$f_{m_0}$} (X6);
         \draw[<-] (X6) to node {$f_{m-1}$} (X7);
          \draw[<-] (X7) to node {$p_{m}$} (X9);
\node (Y1) [below of=Z1] {$Y_1$};
  \node (Y2) [right of=Y1] {$Y_{2}$};
    \node (Y3) [right of=Y2] {$\ldots$};
     \node (Y4) [right of=Y3] {$Y_n$};
      \node (Y5) [right of=Y4] {};
      \node (Y6) [below of=Z9] {$\varprojlim\{Y_\ell,g_\ell\}_{\ell=1}^{\infty}$};
  \draw[<-] (Y1) to node {$g_1$} (Y2);
   \draw[<-] (Y2) to node {$g_2$} (Y3);
     \draw[<-] (Y3) to node {$g_{n-1}$} (Y4);
      \draw[<-] (Y4) to node {$q_n$} (Y6);
      \draw[o-] (Y4) to node {$T_{n,m}(g)$} (X7);
       \draw[<-] (Y6) to node {$g$} (X9);
\end{tikzpicture}
\caption{The diagram from Theorem \ref{Mio:USC}} \label{diagram4}
\end{figure}
\end{theorem}
\begin{proof}
We choose and fix  any $\delta>0$ such that for each set $A\subseteq  \varprojlim\{X_\ell,f_\ell\}_{\ell=1}^{\infty}$, 
$$
\diam(A)<\delta \Longrightarrow \diam(q_n(g(A)))<\varepsilon
$$
holds. Such a $\delta$ does exist since $q_n\circ g$ is uniformly continuous. 
Let $m_0$ be a positive integer such that $\frac{1}{2^{m_0}}<\delta$ and let $m\geq m_0$. 

Let $\mathbf{x}\in \varprojlim\{X_\ell,f_\ell\}_{\ell=1}^{\infty}$ be any point.  Then 
$$
T_{n,m}(g)(p_m(\mathbf x))=q_n(g(p_m^{-1}(p_m(\mathbf x)))),
$$
and, since $\diam(p_m^{-1}(p_m(\mathbf x)))\leq \frac{1}{2^m}\leq \frac{1}{2^{m_0}}<\delta$, it follows that 
$$
\diam(T_{n,m}(g)(p_m(\mathbf x)))<\varepsilon.
$$
Since $\mathbf x\in p_m^{-1}(p_m(\mathbf x))$, it follows that
$$
q_n(g(\mathbf x))\in q_n(g(p_m^{-1}(p_m(\mathbf x)))) = T_{n,m}(g)(p_m(\mathbf x)).
$$
\end{proof}

\begin{theorem}\label{Mio:USC1}
Let $\{X_\ell,f_\ell\}_{\ell=1}^{\infty}$ and $\{Y_\ell,g_\ell\}_{\ell=1}^{\infty}$ be inverse sequences of compact metric spaces and surjective continuous bonding functions, let $g:\varprojlim\{X_\ell,f_\ell\}_{\ell=1}^{\infty}\rightarrow \varprojlim\{Y_\ell,g_\ell\}_{\ell=1}^{\infty}$ be a continuous function,
and let $\varepsilon>0$ and $n$ a positive integer.  Then there is a positive integer $m_0$ such that for each positive integer $m\geq m_0$,  
\begin{enumerate}
\item $\diam (g_{j,n}(T_{n,m}(g)(x)))<\varepsilon $ for any $x\in X_m$,
\item $q_j(g(\mathbf x))\in  g_{j,n}(T_{n,m}(g)(p_m(\mathbf x)))$ for each $\mathbf x\in\varprojlim\{X_\ell,f_\ell\}_{\ell=1}^{\infty}$
\end{enumerate}
for any $j\in\{1,2,3,\ldots,n\}$. [See Figure \ref{diagram5}.]

\begin{figure}[h]
\centering
\noindent\begin{tikzpicture}[node distance=1.5cm, auto]
  \node (X1) {$X_1$};
  \node (X2) [right of=X1] {$X_{2}$};
    \node (X3) [right of=X2] {$X_{3}$};
    \node (X4) [right of=X3] {$\ldots$};
      \node (X5) [right of=X4] {$X_{m_0}$};
    \node (X6) [right of=X5] {$\ldots$};
        \node (X7) [right of=X6] {$X_{m}$};
    \node (X8) [right of=X7] {};
    \node (X9) [right of=X8] {$\varprojlim\{X_\ell,f_\ell\}_{\ell=1}^{\infty}$};
   \node (Z1) [below of=X1] {};
  \node (Z2) [right of=Z1] {};
    \node (Z3) [right of=Z2] {};
    \node (Z4) [right of=Z3] {};
      \node (Z5) [right of=Z4] {};
    \node (Z6) [right of=Z5] {};
        \node (Z7) [right of=Z6] {};
    \node (Z8) [right of=Z7] {};
    \node (Z9) [right of=Z8] {};
  \draw[<-] (X1) to node {$f_1$} (X2);
   \draw[<-] (X2) to node {$f_2$} (X3);
     \draw[<-] (X3) to node {$f_3$} (X4);
   \draw[<-] (X4) to node {$f_{m_0-1}$} (X5);
      \draw[<-] (X5) to node {$f_{m_0}$} (X6);
         \draw[<-] (X6) to node {$f_{m-1}$} (X7);
          \draw[<-] (X7) to node {$p_{m}$} (X9);
\node (Y1) [below of=Z1] {$Y_1$};
  \node (Y2) [right of=Y1] {$Y_{2}$};
    \node (Y3) [right of=Y2] {$\ldots$};
      \node (a) [right of=Y3] {$Y_j$};
       \node (b) [right of=a] {};
     \node (Y4) [right of=b] {$Y_n$};
      \node (Y5) [right of=Y4] {};
      \node (Y6) [below of=Z9] {$\varprojlim\{Y_\ell,g_\ell\}_{\ell=1}^{\infty}$};
  \draw[<-] (Y1) to node {$g_1$} (Y2);
   \draw[<-] (Y2) to node {$g_2$} (Y3);
     \draw[<-] (a) to node {$g_{j,n}$} (Y4);
   \draw[<-] (Y3) to node {$g_{j-1}$} (a);
      \draw[<-] (Y4) to node {$q_n$} (Y6);
      \draw[o-] (Y4) to node {$T_{n,m}(g)$} (X7);
       \draw[<-] (Y6) to node {$g$} (X9);
\end{tikzpicture}
\caption{The diagram from Theorem \ref{Mio:USC1}} \label{diagram5}
\end{figure}
\end{theorem}
\begin{proof}
Let $j\in\{1,2,3,\ldots,n\}$.  Choose any $\delta>0$ such that for each set $A\subseteq  Y_n$, 
$$
\diam(A)<\delta \Longrightarrow \diam(g_{j,n}(A))<\varepsilon
$$
holds. Such a $\delta$ does exist since $g_{j,n}$ is uniformly continuous. 
Let $m_0$ be a positive integer such that for each positive integer $m\geq m_0$,  the following holds
\begin{enumerate}
\item $\diam (T_{n,m}(g)(x))<\delta $ for any $x\in X_m$,
\item $q_n(g(\mathbf x))\in  T_{n,m}(g)(p_m(\mathbf x))$ for each $\mathbf x\in\varprojlim\{X_\ell,f_\ell\}_{\ell=1}^{\infty}$. 
\end{enumerate}
Note that such an $m_0$ does exist by Theorem \ref{Mio:USC}. 
Let $m$ be a positive integer such that $m\geq m_0$.

It follows that
\begin{enumerate}
\item $\diam (g_{j,n}(T_{n,m}(g)(x)))<\varepsilon $ for any $x\in X_m$,
\item $g_{j,n}(q_n(g(\mathbf x)))=q_j(g(\mathbf x))\in  g_{j,n}(T_{n,m}(g)(p_m(\mathbf x)))$ for each $\mathbf x\in\varprojlim\{X_\ell,f_\ell\}_{\ell=1}^{\infty}$.
\end{enumerate}
\end{proof}
In the following theorem a relation is obtained between a continuous function $g:\varprojlim\{X_\ell,f_\ell\}_{\ell=1}^{\infty}\rightarrow \varprojlim\{Y_\ell,g_\ell\}_{\ell=1}^{\infty}$ and a family of $T_{n,m}$-transformations $T_{n_k,m_k}(g)$. 
\begin{theorem}\label{Mio:USC2}
Let $\{X_\ell,f_\ell\}_{\ell=1}^{\infty}$ and $\{Y_\ell,g_\ell\}_{\ell=1}^{\infty}$ be inverse sequences of compact metric spaces and surjective continuous bonding functions, and let $g:\varprojlim\{X_\ell,f_\ell\}_{\ell=1}^{\infty}\rightarrow \varprojlim\{Y_\ell,g_\ell\}_{\ell=1}^{\infty}$ be a continuous function. 
Then there are  sequences $(n_\ell)$ and $(m_\ell)$ of positive integers such that 
\begin{enumerate}
\item \label{unounouno}$m_{k+1}\geq m_{k}$ and $n_{k+1} > n_{k}$ for each positive integer $k$,
\item \label{unouno}$m_k\geq n_k$ for each positive integer $k$,
\item \label{uno} for each positive integer $k$ and for each $j\in\{1,2,3,\ldots ,n_k\}$, 
$$
g_{j,n_r}(T_{n_r,m_r}(g)(x))\subseteq g_{j,n_k}(T_{n_k,m_k}(g)(f_{m_k,m_r}(x)))
$$
holds for each positive integer $r > k$  and  for each $x\in X_{m_r}$ (see Figure \ref{diagram6}),

\begin{figure}[h]
\centering
\begin{tikzpicture}[node distance=1.5cm, auto]
  \node (X1) {$X_{m_k}$};
  \node (X2) [right of=X1] {};
  \node (X3) [right of=X2] {$X_{m_r}$};
  \draw[<-] (X1) to node {$f_{m_k,m_r}$} (X3);
  \node (Z) [below of=X1] {};
  \node (Y3) [below of=Z] {$Y_{n_k}$};
\node (Y2) [left of=Y3] {};
 \node (Y1) [left of=Y2] {$Y_{j}$};
 \node (Y4) [right of=Y3] {};
  \node (Y5) [right of=Y4] {$Y_{n_r}$};
  \draw[<-] (Y1) to node {$g_{j,n_k}$} (Y3);
   \draw[<-] (Y3) to node {$g_{n_k,n_r}$} (Y5);
               \draw[o-] (Y3) to node {$T_{n_k,m_k}(g)$} (X1);
                \draw[o-] (Y5) to node {$T_{n_r,m_r}(g)$} (X3);
\end{tikzpicture}
\caption{The diagram from item \ref{uno}} \label{diagram6}
\end{figure}
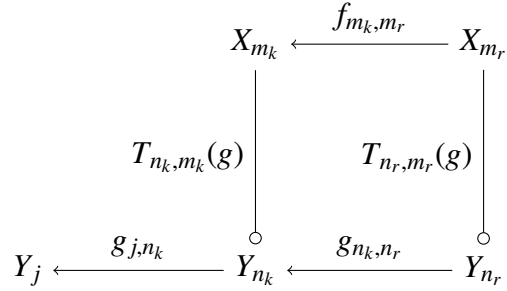

\item \label{due} for each positive integer $k$ and for each $j\in\{1,2,3,\ldots ,n_k\}$, 
$$
q_{j}(g(\mathbf x))\in g_{j,n_k}(T_{n_k,m_k}(g)(p_{m_k}(\mathbf x)))
$$
holds for each $\mathbf x\in \varprojlim\{X_\ell,f_\ell\}_{\ell=1}^{\infty}$ (see Figure \ref{diagram7}),

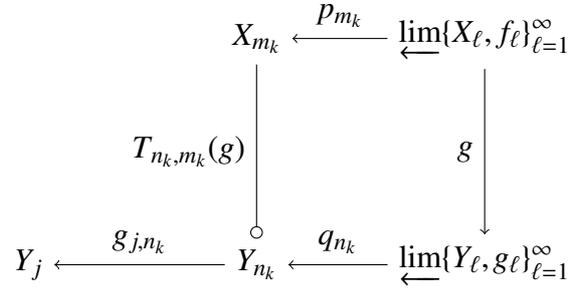
\begin{figure}[h]
\centering
\begin{tikzpicture}[node distance=1.5cm, auto]
  \node (X1) {$X_{m_k}$};
  \node (X2) [right of=X1] {};
  \node (X3) [right of=X2] {$\varprojlim\{X_\ell,f_\ell\}_{\ell=1}^{\infty}$};
  \draw[<-] (X1) to node {$p_{m_k}$} (X3);
   \node (Z) [below of=X1] {};
    \node (Y3) [below of=Z] {$Y_{n_k}$};
      \node (Y2) [left of=Y3] {};
        \node (Y1) [left of=Y2] {$Y_{j}$};
          \node (Y4) [right of=Y3] {};
            \node (Y5) [right of=Y4] {$\varprojlim\{Y_\ell,g_\ell\}_{\ell=1}^{\infty}$};
             \draw[<-] (Y1) to node {$g_{j,n_k}$} (Y3);
              \draw[<-] (Y3) to node {$q_{n_k}$} (Y5);
               \draw[o-] (Y3) to node {$T_{n_k,m_k}(g)$} (X1);
                \draw[<-] (Y5) to node {$g$} (X3);
\end{tikzpicture}
\caption{The diagram from item \ref{due}} \label{diagram7}
\end{figure}

\item \label{tre} for each positive integer $j$, 
$$
\lim_{k\to\infty}\diam (g_{j,n_k}(T_{n_k,m_k}(g)(p_{m_k}(\mathbf x))))=0
$$
for any $\mathbf x\in \varprojlim\{X_\ell,f_\ell\}_{\ell=1}^{\infty}$.
\end{enumerate}
\end{theorem}
\begin{proof}
We are constructing a diagram as presented on Figure \ref{diagram8} such that \ref{unounouno}, \ref{unouno}, \ref{uno}, \ref{due} and \ref{tre} are satisfied. 

\begin{figure}[h]
\centering
\noindent \begin{tikzpicture}[node distance=1.5cm, auto]
  \node (X1) {$X_{m_1}$};
  \node (X2) [right of=X1] {};
   \node (X3) [right of=X2] {$X_{m_2}$};
    \node (X4) [right of=X3] {};
      \node (X5) [right of=X4] {$X_{m_3}$};
    \node (X6) [right of=X5] {};
        \node (X7) [right of=X6] {$X_{m_4}$};
    \node (X8) [right of=X7] {$\ldots$};
    \node (X9) [right of=X8] {$\varprojlim\{X_\ell,f_\ell\}_{\ell=1}^{\infty}$};
   \node (Z1) [below of=X1] {};
  \node (Z2) [right of=Z1] {};
    \node (Z3) [right of=Z2] {};
    \node (Z4) [right of=Z3] {};
      \node (Z5) [right of=Z4] {};
    \node (Z6) [right of=Z5] {};
        \node (Z7) [right of=Z6] {};
    \node (Z8) [right of=Z7] {};
    \node (Z9) [right of=Z8] {};
  \draw[<-] (X1) to node {$f_{m_1,m_2}$} (X3);
   \draw[<-] (X3) to node {$f_{m_2,m_3}$} (X5);
      \draw[<-] (X5) to node {$f_{m_3,m_4}$} (X7);
          \draw[<-] (X7) to node {$f_{m_4}$} (X8);
\node (Y1) [below of=Z1] {$Y_{n_1}$};
  \node (Y2) [right of=Y1] {};
    \node (Y3) [right of=Y2] {$Y_{n_2}$};
      \node (a) [right of=Y3] {};
       \node (b) [right of=a] {$Y_{n_3}$};
     \node (Y4) [right of=b] {};
      \node (Y5) [right of=Y4] {$Y_{n_4}$};
         \node (Y6) [right of=Y5] {$\ldots$};
      \node (Y7) [below of=Z9] {$\varprojlim\{Y_\ell,g_\ell\}_{\ell=1}^{\infty}$};
  \draw[<-] (Y1) to node {$g_{n_1,n_2}$} (Y3);
   \draw[<-] (Y3) to node {$g_{n_2,n_3}$} (b);
     \draw[<-] (b) to node {$g_{n_3,n_4}$} (Y5);
   \draw[<-] (Y5) to node {$g_{n_4}$} (Y6);
      \draw[o-] (Y1) to node {$T_{n_1,m_1}(g)$} (X1);
            \draw[o-] (Y3) to node {$T_{n_2,m_2}(g)$} (X3);
                  \draw[o-] (b) to node {$T_{n_3,m_3}(g)$} (X5);
                        \draw[o-] (Y5) to node {$T_{n_4,m_4}(g)$} (X7);
       \draw[<-] (Y7) to node {$g$} (X9);
\end{tikzpicture}
\caption{Constructing a diagram such that \ref{unounouno}, \ref{unouno}, \ref{uno}, \ref{due} and \ref{tre} are satisfied} \label{diagram8}
\end{figure}
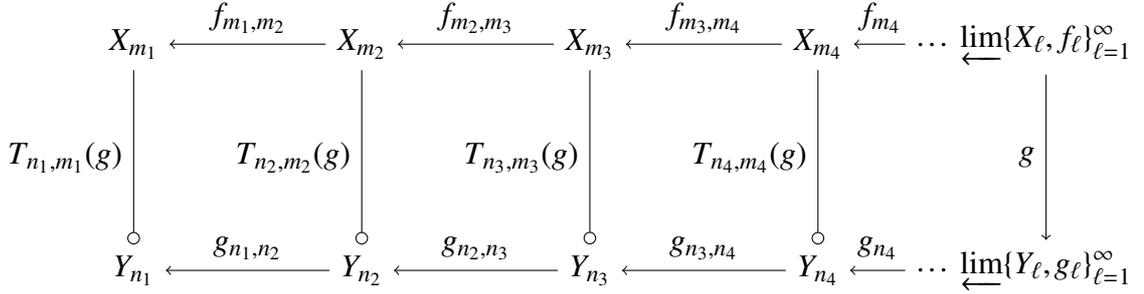

First, let $(\varepsilon_\ell)$ be a sequence of positive real numbers such that $\displaystyle \lim_{\ell\to \infty}\varepsilon_\ell=0$. 

Let $n_1$ be any positive integer.  By Theorem \ref{Mio:USC1}, there is a positive integer  $m_1$ such that for each positive integer $m\geq m_1$ and for any $j\in\{1,2,3,\ldots,n_1\}$,  
\begin{enumerate}
\item $\diam (g_{j,n_1}(T_{n_1,m}(g)(x)))<\varepsilon_1 $ for any $x\in X_{m}$,
\item $q_j(g(\mathbf x))\in  g_{j,n_1}(T_{n_1,m}(g)(p_m(\mathbf x)))$ for each $\mathbf x\in\varprojlim\{X_\ell,f_\ell\}_{\ell=1}^{\infty}$.  
\end{enumerate}
 Choose and fix such a positive integer $m_1\geq n_1$.

Obviously,  \ref{due}  for $k=1$ and $j\in \{1,2,3,\ldots ,n_1\}$ holds, i.e., we have that 
$$
q_{j}(g(\mathbf x))\in g_{j,n_1}(T_{n_1,m_1}(g)(p_{m_1}(\mathbf x)))
$$
holds for each $\mathbf x\in \varprojlim\{X_\ell,f_\ell\}_{\ell=1}^{\infty}$.

Suppose next, that for a positive integer $k_0$, we have already constructed $k_0$-tuples  $(m_1,m_2,m_3,\ldots ,m_{k_0})$ and $(n_1,n_2,n_3,\ldots ,n_{k_0})$ of positive integers such that 
\begin{enumerate}
\item $m_{k+1}\geq m_{k}$ and $n_{k+1} > n_{k}$ for each  $k\in \{1,2,3,\ldots ,k_0-1\}$,
\item $m_k\geq n_k$ for each for each  $k\in \{1,2,3,\ldots ,k_0\}$,
\item  for each  $k\in \{1,2,3,\ldots ,k_0-1\}$ and for each $j\in\{1,2,3,\ldots ,n_k\}$, 
$$
g_{j,n_r}(T_{n_r,m_r}(g)(x))\subseteq g_{j,n_k}(T_{n_k,m_k}(g)(f_{m_k,m_r}(x)))
$$
holds for each positive integer $r \in \{k+1,k+2,k+3,\ldots ,k_0\}$  and  for each $x\in X_{m_r}$,
\item  for each  $k\in \{1,2,3,\ldots ,k_0\}$ and for each $j\in\{1,2,3,\ldots ,n_k\}$,
$$
g_{j,n_k}(q_{n_k}(g(\mathbf x)))\in g_{j,n_k}(T_{n_k,m_k}(g)(p_{m_k}(\mathbf x)))
$$
holds for each $\mathbf x\in \varprojlim\{X_\ell,f_\ell\}_{\ell=1}^{\infty}$,
\item for each  $k\in \{1,2,3,\ldots ,k_0\}$ and for each $j\in\{1,2,3,\ldots ,n_k\}$, 
$$
\diam (g_{j,n_k}(T_{n_k,m_k}(g)(x)))<\varepsilon_k 
$$
for any $x\in X_{m_k}$.
\end{enumerate}
Let $n_{{k_0}+1}>n_{k_0}$ be any positive integer.  By Theorem \ref{Mio:USC1}, there is a positive integer  $m_{{k_0}+1}$ such that for each positive integer $m\geq m_{{k_0}+1}$ and for any $j\in\{1,2,3,\ldots,n_{k_0+1}\}$,  
\begin{enumerate}
\item $\diam (g_{j,n_{{k_0}+1}}(T_{{n_{{k_0}+1}},m}(g)(x)))<\varepsilon_{{k_0}+1} $ for any $x\in X_{m}$,
\item $q_j(g(\mathbf x))\in  g_{j,n_{{k_0}+1}}(T_{{n_{{k_0}+1}},m}(g)(p_m(\mathbf x)))$ for each $\mathbf x\in\varprojlim\{X_\ell,f_\ell\}_{\ell=1}^{\infty}$.
\end{enumerate}
 Choose and fix such a positive integer $m_{{k_0}+1}\geq n_{{k_0}+1}$.

First, we prove \ref{uno} for $r={k_0}+1$, $k\in \{1,2,3,\ldots ,k_0\}$ and $j \in \{1,2,3,\ldots ,n_k\}$, i.e.,  that 
$$
g_{j,n_{{k_0}+1}}(T_{{n_{{k_0}+1}},m_{{k_0}+1}}(g)(x))\subseteq g_{j,n_k}(T_{n_k,m_k}(g)(f_{m_k,m_{{k_0}+1}}(x)))
$$
  for each $x\in X_{m_{{k_0}+1}}$, see Figure \ref{diagram9}.

\begin{figure}[h]
\centering
\begin{tikzpicture}[node distance=1.5cm, auto]
  \node (X1) {$X_{m_k}$};
  \node (X2) [right of=X1] {};
  \node (X3) [right of=X2] {$X_{m_{{k_0}+1}}$};
  \draw[<-] (X1) to node {$f_{m_k,m_{{k_0}+1}}$} (X3);
   \node (Z) [below of=X1] {};
    \node (Y3) [below of=Z] {$Y_{n_k}$};
      \node (Y2) [left of=Y3] {};
        \node (Y1) [left of=Y2] {$Y_{j}$};
          \node (Y4) [right of=Y3] {};
            \node (Y5) [right of=Y4] {$Y_{n_{{k_0}+1}}$};
             \draw[<-] (Y1) to node {$g_{j,n_k}$} (Y3);
              \draw[<-] (Y3) to node {$g_{n_k,n_{{k_0}+1}}$} (Y5);
               \draw[o-] (Y3) to node {$T_{{n_{k},m_{k}}}(g)$} (X1);
                \draw[o-] (Y5) to node {$T_{{n_{{k_0}+1}},m_{{k_0}+1}}(g)$} (X3);
\end{tikzpicture}
\caption{The diagram from item \ref{uno} for $r={k_0}+1$, $k\in \{1,2,3,\ldots ,k_0\}$ and $j \in \{1,2,3,\ldots ,n_k\}$} \label{diagram9}
\end{figure}

Let $x\in X_{m_{{k_0}+1}}$ be any point.  First, we prove that 
$$
p_{m_{{k_0}+1}}^{-1}(x)\subseteq p_{m_k}^{-1}(f_{m_k,m_{{k_0}+1}}(x)).
$$
Let $\mathbf x\in p_{m_{{k_0}+1}}^{-1}(x)$. Then $p_{m_{{k_0}+1}}(\mathbf x)=x$ and therefore
$$
p_{m_k}(\mathbf x)=f_{m_k,m_{{k_0}+1}}(p_{m_{{k_0}+1}}(\mathbf x))=f_{m_k,m_{{k_0}+1}}(x),
$$
proving that $\mathbf x\in p_{m_k}^{-1}(f_{m_k,m_{{k_0}+1}}(x))$ and $p_{m_{{k_0}+1}}^{-1}(x)\subseteq p_{m_k}^{-1}(f_{m_k,m_{{k_0}+1}}(x))$ follows.
It follows from 
$$
g_{j,n_{{k_0}+1}}(T_{{n_{{k_0}+1}},m_{{k_0}+1}}(g)(x))=g_{j,n_k}(g_{n_k,n_{{k_0}+1}}(q_{n_{{k_0}+1}}( g( p_{m_{{k_0}+1}}^{-1}(x)))))
$$
$$
=g_{j,n_k}(q_{n_k}( g( p_{m_{{k_0}+1}}^{-1}(x))))
$$
and 
$$
g_{j,n_k}(T_{n_{k},m_{k}}(g)(f_{m_k,m_{{k_0}+1}}(x)))=g_{j,n_k}(q_{n_k}(g( p_{m_{k}}^{-1}(f_{m_{k},m_{{k_0}+1}}(x)))))
$$
that 
$$
g_{j,n_{{k_0}+1}}(T_{{n_{{k_0}+1}},m_{{k_0}+1}}(g)(x))\subseteq g_{j,n_k}(T_{n_{k},m_{k}}(g)(f_{m_k,m_{{k_0}+1}}(x))).
$$
Note that \ref{due} for $k=k_0+1$ and $j\in \{1,2,3,\ldots ,n_k\}$ is also satisfied, i.e., that  
$$
q_{j}(g(\mathbf x))\in g_{j,n_{k_0+1}}(T_{{n_{{k_0}+1}},m_{{k_0}+1}}(g)(p_{m_{k_0+1}}(\mathbf x)))
$$
is satisfied for each $\mathbf x\in \varprojlim\{X_\ell,f_\ell\}_{\ell=1}^{\infty}$ (see Figure \ref{diagram10}).

\begin{figure}[h]
\centering
\begin{tikzpicture}[node distance=1.5cm, auto]
  \node (X1) {$X_{m_{k_0+1}}$};
  \node (X2) [right of=X1] {};
  \node (X3) [right of=X2] {$\varprojlim\{X_\ell,f_\ell\}_{\ell=1}^{\infty}$};
  \draw[<-] (X1) to node {$p_{m_{k_0+1}}$} (X3);
   \node (Z) [below of=X1] {};
    \node (Y3) [below of=Z] {$Y_{n_{k_0+1}}$};
      \node (Y2) [left of=Y3] {};
        \node (Y1) [left of=Y2] {$Y_{j}$};
          \node (Y4) [right of=Y3] {};
            \node (Y5) [right of=Y4] {$\varprojlim\{Y_\ell,g_\ell\}_{\ell=1}^{\infty}$};
             \draw[<-] (Y1) to node {$g_{j,n_{k_0+1}}$} (Y3);
              \draw[<-] (Y3) to node {$q_{n_{k_0+1}}$} (Y5);
               \draw[o-] (Y3) to node {$T_{{n_{{k_0}+1}},m_{{k_0}+1}}(g)$} (X1);
                \draw[<-] (Y5) to node {$g$} (X3);
\end{tikzpicture}
\caption{The diagram from item \ref{due} for $k=k_0+1$ and $j\in \{1,2,3,\ldots ,n_k\}$} \label{diagram10}
\end{figure}

Finally, note that for any $j\in\{1,2,3,\ldots,n_{k_0+1}\}$,
 $$
 \diam (g_{j,n_{k_0+1}}(T_{{n_{{k_0}+1}},m_{{k_0}+1}}(g)(x)))<\varepsilon_{k_0+1} 
 $$
 for any $x\in X_{m_{k_0+1}}$. 
 
 Inductively, we have constructed sequences $(m_\ell)$ and $(n_\ell)$ of positive integers such that
 \begin{enumerate}
\item $m_{k+1}\geq m_{k}$ and $n_{k+1} > n_{k}$ for each positive integer $k$,
\item $m_k\geq n_k$ for each positive integer $k$,
\item for each positive integer $k$ and for each $j\in\{1,2,3,\ldots ,n_k\}$, 
$$
g_{j,n_r}(T_{n_r,m_r}(g)(x))\subseteq g_{j,n_k}(T_{n_k,m_k}(g)(f_{m_k,m_r}(x)))
$$
holds for each positive integer $r > k$  and  for each $x\in X_{m_r}$,
\item for each positive integer $k$ and for each $j\in\{1,2,3,\ldots ,n_k\}$, 
$$
q_{j}(g(\mathbf x))\in g_{j,n_k}(T_{n_k,m_k}(g)(p_{m_k}(\mathbf x)))
$$
holds for each $\mathbf x\in \varprojlim\{X_\ell,f_\ell\}_{\ell=1}^{\infty}$, and
\item  for each positive integer $k$ and for each $j\in\{1,2,3,\ldots ,n_k\}$, 
$$
\diam (g_{j,n_k}(T_{n_k,m_k}(g)(x)))<\varepsilon_k 
$$
for any $x\in X_{m_k}$.
\end{enumerate}
Since $\displaystyle \lim_{k\to \infty}\varepsilon_k=0$, it follows that 
for each positive integer $j$, 
$$
\lim_{k\to\infty}\diam (g_{j,n_k}(T_{n_k,m_k}(g)(p_{m_k}(\mathbf x))))=0
$$
for any $\mathbf x\in \varprojlim\{X_\ell,f_\ell\}_{\ell=1}^{\infty}$. 
 \end{proof}
Among other things, the following theorem gives sufficient conditions for the existence of a continuous surjection  $\varprojlim\{X_\ell,f_\ell\}_{\ell=1}^{\infty}\rightarrow \varprojlim\{Y_\ell,g_\ell\}_{\ell=1}^{\infty}$. 
\begin{theorem}\label{jejhata}
Let $\{X_\ell,f_\ell\}_{\ell=1}^{\infty}$ and $\{Y_\ell,g_\ell\}_{\ell=1}^{\infty}$ be inverse sequences of compact metric spaces and surjective continuous bonding functions, 
 let $(n_\ell)$ and $(m_\ell)$ be any sequences of positive integers 
such that 
\begin{enumerate}
\item $m_{k+1}\geq m_{k}$ and $n_{k+1} > n_{k}$ for each positive integer $k$,
\item $m_k\geq n_k$ for each positive integer $k$,
\end{enumerate}
and let $(H_\ell)$ be any sequence of upper semicontinuous functions $H_\ell:X_{m_\ell}\multimap Y_{n_\ell}$ such that 
\begin{enumerate}
\item[(a)] \label{uno1} for each positive integer $k$ and for each $j\in\{1,2,3,\ldots ,n_k\}$, 
$$
g_{j,n_r}(H_r(x))\subseteq g_{j,n_k}(H_k(f_{m_k,m_r}(x)))
$$
holds for each positive integer $r > k$  and  for each $x\in X_{m_r}$ (see Figure \ref{diagram11}),

\begin{figure}[h]
\centering
\begin{tikzpicture}[node distance=1.5cm, auto]
  \node (X1) {$X_{m_k}$};
  \node (X2) [right of=X1] {};
  \node (X3) [right of=X2] {$X_{m_r}$};
  \draw[<-] (X1) to node {$f_{m_k,m_r}$} (X3);
   \node (Z) [below of=X1] {};
    \node (Y3) [below of=Z] {$Y_{n_k}$};
      \node (Y2) [left of=Y3] {};
        \node (Y1) [left of=Y2] {$Y_{j}$};
          \node (Y4) [right of=Y3] {};
            \node (Y5) [right of=Y4] {$Y_{n_r}$};
             \draw[<-] (Y1) to node {$g_{j,n_k}$} (Y3);
              \draw[<-] (Y3) to node {$g_{n_k,n_r}$} (Y5);
               \draw[o-] (Y3) to node {$H_k$} (X1);
                \draw[o-] (Y5) to node {$H_r$} (X3);
\end{tikzpicture}
\caption{Fhe diagram from item $(a)$} \label{diagram11}
\end{figure}

\item[(b)] \label{tre3} for each positive integer $j$, 
$$
\lim_{k\to\infty}\diam (g_{j,n_k}(H_{k}(p_{m_k}(\mathbf x))))=0
$$
for any $\mathbf x\in \varprojlim\{X_\ell,f_\ell\}_{\ell=1}^{\infty}$. 
\end{enumerate}
Then the following holds true. 
\begin{enumerate}
\item\label{jej1}
There is a continuous function $g:\varprojlim\{X_\ell,f_\ell\}_{\ell=1}^{\infty}\rightarrow \varprojlim\{Y_\ell,g_\ell\}_{\ell=1}^{\infty}$ such that for each positive integer $k$ and for each $j\in \{1,2,3,\ldots, n_k\}$,
$$
q_j(g(\mathbf x))\in g_{j,n_k}( H_k (p_{m_k}(\mathbf x)))
$$
for each $\mathbf x\in \varprojlim\{X_\ell,f_\ell\}_{\ell=1}^{\infty}$.
\item \label{jej2} If the graph of each $H_k$ is surjective, then  there is a continuous surjection $g:\varprojlim\{X_\ell,f_\ell\}_{\ell=1}^{\infty}\rightarrow \varprojlim\{Y_\ell,g_\ell\}_{\ell=1}^{\infty}$ such that for each positive integer $k$ and for each $j\in \{1,2,3,\ldots, n_k\}$,
$$
q_j(g(\mathbf x))\in g_{j,n_k}( H_k (p_{m_k}(\mathbf x)))
$$
for each $\mathbf x\in \varprojlim\{X_\ell,f_\ell\}_{\ell=1}^{\infty}$.
\item \label{jej3} If $f,g:\varprojlim\{X_\ell,f_\ell\}_{\ell=1}^{\infty}\rightarrow \varprojlim\{Y_\ell,g_\ell\}_{\ell=1}^{\infty}$ are continuous functions such that for any positive integer $k$ and for any $j\in \{1,2,3,\ldots, n_k\}$, 
$$
q_{j}(g(\mathbf x))\in g_{j,n_k}(H_k(p_{m_k}(\mathbf x))) \textup{ and }  q_{j}(f(\mathbf x))\in g_{j,n_k}(H_k(p_{m_k}(\mathbf x)))
$$
holds for each $\mathbf x\in \varprojlim\{X_\ell,f_\ell\}_{\ell=1}^{\infty}$, then $f=g$ (see Figure \ref{diagram12}).

\begin{figure}[h]
\centering
\begin{tikzpicture}[node distance=1.5cm, auto]
  \node (X1) {$X_{m_k}$};
  \node (X2) [right of=X1] {};
  \node (X3) [right of=X2] {$\varprojlim\{X_\ell,f_\ell\}_{\ell=1}^{\infty}$};
  \draw[<-] (X1) to node {$p_{m_k}$} (X3);
   \node (Z) [below of=X1] {};
    \node (Y3) [below of=Z] {$Y_{n_k}$};
      \node (Y2) [left of=Y3] {};
        \node (Y1) [left of=Y2] {$Y_{j}$};
          \node (Y4) [right of=Y3] {};
            \node (Y5) [right of=Y4] {$\varprojlim\{Y_\ell,g_\ell\}_{\ell=1}^{\infty}$};
             \draw[<-] (Y1) to node {$g_{j,n_k}$} (Y3);
              \draw[<-] (Y3) to node {$q_{n_k}$} (Y5);
               \draw[o-] (Y3) to node {$H_k$} (X1);
               \draw[<-] (Y5) to node {$g,f$} (X3);
\end{tikzpicture}
\caption{The diagram from item \ref{jej3}} \label{diagram12}
\end{figure}
\end{enumerate}
\end{theorem}
\begin{proof}
We first prove \ref{jej1}. 
For each $\mathbf x=(x_1,x_2,x_3,\ldots)\in \varprojlim\{X_\ell,f_\ell\}_{\ell=1}^{\infty}$ and for each positive integer $k$, let 
$$
G_k(\mathbf x)=\left(\{(g_{1,n_k}(y),g_{2,n_k}(y),g_{3,n_k}(y),\ldots, g_{n_k-1,n_k}(y),y) \ | \  y\in H_k(x_{m_k})\}\times \prod_{\ell=n_k+1}^{\infty}Y_\ell\right)\cap
$$
$$
\varprojlim\{Y_\ell,g_\ell\}_{\ell=1}^{\infty}.
$$
We prove that $G_k:\varprojlim\{X_\ell,f_\ell\}_{\ell=1}^{\infty}\multimap \varprojlim\{Y_\ell,g_\ell\}_{\ell=1}^{\infty}$ is a well-defined upper semicontinuous set-valued function for each positive integer $k$. Let $k$ be a positive integer and let $\mathbf x\in \varprojlim\{X_\ell,f_\ell\}_{\ell=1}^{\infty}$ be any point.  Obviously, $G_k(\mathbf x)$ is non-empty and closed in $\varprojlim\{Y_\ell,g_\ell\}_{\ell=1}^{\infty}$, therefore $G_k:\varprojlim\{X_\ell,f_\ell\}_{\ell=1}^{\infty}\multimap \varprojlim\{Y_\ell,g_\ell\}_{\ell=1}^{\infty}$ is a well-defined set-valued function. To see that it is an  upper semicontinuous function, we prove that its graph $\Gamma(G_k)$ is closed in $\varprojlim\{X_\ell,f_\ell\}_{\ell=1}^{\infty}\times \varprojlim\{Y_\ell,g_\ell\}_{\ell=1}^{\infty}$. For each positive integer $\ell$, let
$$
(\mathbf x^\ell,\mathbf y^\ell)=((x_1^\ell, x_2^\ell, x_3^\ell, \ldots ),(y_1^\ell, y_2^\ell, y_3^\ell, \ldots ))\in \Gamma(G_k)
$$
and let 
$$
(\mathbf x^0,\mathbf y^0)=((x_1^0, x_2^0, x_3^0, \ldots ),(y_1^0, y_2^0, y_3^0, \ldots )
)\in \varprojlim\{X_\ell,f_\ell\}_{\ell=1}^{\infty}\times \varprojlim\{Y_\ell,g_\ell\}_{\ell=1}^{\infty}
$$
be such that 
$$
\lim_{\ell \to \infty }(\mathbf x^\ell,\mathbf y^\ell)=(\mathbf x^0,\mathbf y^0).
$$
 Note that for each positive integer $\ell$, 
 $$
 y_{n_k}^\ell\in H_k(x_{m_k}^\ell).
$$
Since $H_k$ is an upper semicontinuous set-valued function and since $\displaystyle \lim_{\ell\to \infty}(x_{m_k}^\ell,y_{n_k}^\ell)=(x_{m_k}^0,y_{n_k}^0) $, it follows that $y_{n_k}^0\in H_k(x_{m_k}^0)$. Also, note that
$$
(y_1^0, y_2^0, y_3^0, \ldots ,y_{n_k}^0)=(g_{1,n_k}(y_{n_k}^0), g_{2,n_k}(y_{n_k}^0), g_{3,n_k}(y_{n_k}^0), \ldots , g_{n_k-1,n_k}(y_{n_k}^0),y_{n_k}^0).
$$
Therefore, $(\mathbf x^0,\mathbf y^0)\in \Gamma(G_k)$ and this proves that $G_k$ is upper semicontinuous. 

Note also that for each positive integer $j\leq n_k$,
$$
q_j(G_k(\mathbf x))=g_{j,n_k}(H_k(p_{m_k}(\mathbf x)))
$$
for each $\mathbf x\in \varprojlim\{X_\ell,f_\ell\}_{\ell=1}^{\infty}$ and, therefore, 
$$
\lim_{k\to \infty} \diam(q_j(G_k(\mathbf x)))=\lim_{k\to \infty}\diam (g_{j,n_k}(H_k(p_{m_k}(\mathbf x)))) = 0
$$
follows for each $\mathbf x\in \varprojlim\{X_\ell,f_\ell\}_{\ell=1}^{\infty}$. 
Therefore, for each positive integer $k$,
$$
\lim_{k\to \infty} \diam (G_k(\mathbf x))=0
$$
for each $\mathbf x\in \varprojlim\{X_\ell,f_\ell\}_{\ell=1}^{\infty}$. 
Note also that for each positive integer $k$,
$$
G_{k+1}(\mathbf x)\subseteq G_k(\mathbf x)
$$
for each $\mathbf x\in \varprojlim\{X_\ell,f_\ell\}_{\ell=1}^{\infty}$.
Therefore, for each $\mathbf x\in \varprojlim\{X_\ell,f_\ell\}_{\ell=1}^{\infty}$, there is a point $\mathbf y(\mathbf x)\in \varprojlim\{Y_\ell,g_\ell\}_{\ell=1}^{\infty}$
such that 
$$
\bigcap_{k=1}^{\infty} G_k(\mathbf x)=\{\mathbf y(\mathbf x)\}.
$$
We define
$$
g(\mathbf x)=\mathbf y(\mathbf x)
$$
for each $\mathbf x\in \varprojlim\{X_\ell,f_\ell\}_{\ell=1}^{\infty}$. It follows from General mapping theorem \cite[Theorem 7.4, page 105]{nadler} that $g$ is a well-defined continuous function.  It follows from 
$$
q_j(g(\mathbf x))=g_{j,n_k}(q_{n_k}(g(\mathbf x)))
$$
and 
$$
q_{n_k}(g(\mathbf x))\in H_k(p_{m_k}(\mathbf x))
$$
that 
$$
q_j(g(\mathbf x))\in g_{j,n_k}( H_k (p_{m_k}(\mathbf x)))
$$
for each $\mathbf x\in \varprojlim\{X_\ell,f_\ell\}_{\ell=1}^{\infty}$, for each positive integer $k$ and for each $j\in \{1,2,3,\ldots, n_k\}$.

If, in addition, the graph of each $H_k$ is surjective, then by General mapping theorem, $g$ is also surjective, which proves \ref{jej2}. 

To prove \ref{jej3}, let $f$ and $g$ be any functions  $\varprojlim\{X_\ell,f_\ell\}_{\ell=1}^{\infty}\rightarrow \varprojlim\{Y_\ell,g_\ell\}_{\ell=1}^{\infty}$ such that for any positive integer $k$ and for any $j\in \{1,2,3,\ldots, n_k\}$, 
$$
q_{j}(g(\mathbf x))\in g_{j,n_k}(H_k(p_{m_k}(\mathbf x))) \textup{ and }  q_{j}(f(\mathbf x))\in g_{j,n_k}(H_k(p_{m_k}(\mathbf x)))
$$
holds for each $\mathbf x\in \varprojlim\{X_\ell,f_\ell\}_{\ell=1}^{\infty}$ (see Figure \ref{diagram12}).

Suppose that $f\neq g$.  Then there is a positive integer $j$ and a point $\mathbf x\in \varprojlim\{X_\ell,f_\ell\}_{\ell=1}^{\infty}$ such that 
$$
q_{j}(f(\mathbf x))\neq q_{j}(g(\mathbf x)).
$$
Let $\delta =d_{Y_{j}}(q_{j}(f(\mathbf x)), q_{j}(g(\mathbf x)))$. Then, since $q_{j}(f(\mathbf x)), q_{j}(g(\mathbf x))\in g_{j,n_k}(H_k(p_{m_k}(\mathbf x)))$, it follows that $\diam(g_{j,n_k}(H_k(p_{m_k}(\mathbf x))))\geq \delta$ for each positive integer $k$ such that $n_k>j$. Therefore, 
$$
\lim_{k\to\infty}\diam (g_{j,n_k}(H_{k}(p_{m_k}(\mathbf x))))\neq 0,
$$
which is a contradiction. 
\end{proof}
\begin{observation}\label{tupki}
Assume that all the assumptions from Theorem \ref{jejhata} are satisfied. It follows from \ref{jej3} of Theorem \ref{jejhata} and from the proof of \ref{jej1} of Theorem \ref{jejhata} that if 
$$
g:\varprojlim\{X_\ell,f_\ell\}_{\ell=1}^{\infty}\rightarrow \varprojlim\{Y_\ell,g_\ell\}_{\ell=1}^{\infty}
$$
is a continuous function such that for any positive integer $k$ and for any $j\in \{1,2,3,\ldots, n_k\}$ and for each $\mathbf x\in \varprojlim\{X_\ell,f_\ell\}_{\ell=1}^{\infty}$, it holds that
$$
q_{j}(g(\mathbf x))\in g_{j,n_k}(H_k(p_{m_k}(\mathbf x))), 
$$
  then $g(\mathbf x)$ is the single point in $\bigcap_{k=1}^{\infty} G_k(\mathbf x)$, where 
$$
G_k(\mathbf x)=\left(\{(g_{1,n_k}(y),g_{2,n_k}(y),g_{3,n_k}(y),\ldots, g_{n_k-1,n_k}(y),y) \ | \  y\in H_k(x_{m_k})\}\times \prod_{\ell=n_k+1}^{\infty}Y_\ell\right)\cap
$$
$$
\varprojlim\{Y_\ell,g_\ell\}_{\ell=1}^{\infty}
$$
for each $\mathbf x=(x_1,x_2,x_3,\ldots)\in \varprojlim\{X_\ell,f_\ell\}_{\ell=1}^{\infty}$ and for each positive integer $k$. 

Let $\mathbf x=(x_1,x_2,x_3,\ldots)\in \varprojlim\{X_\ell,f_\ell\}_{\ell=1}^{\infty}$ be any point  and let $j$ be any positive integer. Then (taking into account the definition of each $G_k$),
$$
q_{j}(g(\mathbf x))\in q_{j}\left(\bigcap_{k=1}^{\infty} G_k(\mathbf x)\right)\subseteq \bigcap_{k=1}^{\infty}q_{j}(G_k(\mathbf x))=\bigcap_{k=1}^{\infty}g_{j,n_{k}}(H_{k}(p_{m_{k}}(\mathbf x))).
$$
Since 
$$
\lim_{k\to\infty}\diam (g_{j,n_{k}}(H_{k}(p_{m_{k}}(\mathbf x))))=0
$$
it follows that 
$$
\diam \left(\bigcap_{k=1}^{\infty}g_{j,n_{k}}(H_{k}(p_{m_{k}}(\mathbf x))\right)=0.
$$
Therefore, 
$$
\bigcap_{k=1}^{\infty}g_{j,n_{k}}(H_{k}(p_{m_{k}}(\mathbf x)))=\{q_{j}(g(\mathbf x))\}.
$$
\end{observation}
The following corollary is our first mapping theorem for inverse limits. 
\begin{corollary}\label{mappingTHM}
Let $\{X_\ell,f_\ell\}_{\ell=1}^{\infty}$ and $\{Y_\ell,g_\ell\}_{\ell=1}^{\infty}$ be inverse sequences of compact metric spaces and surjective continuous bonding functions. The following statements are equivalent.
\begin{enumerate}
    \item\label{1ena1} There is a continuous surjection  from $\varprojlim\{X_\ell,f_\ell\}_{\ell=1}^{\infty}$ to $\varprojlim\{Y_\ell,g_\ell\}_{\ell=1}^{\infty}$. 
    \item \label{2ena2} There are sequences $(n_\ell)$ and $(m_\ell)$ of positive integers 
such that 
\begin{enumerate}
\item $m_{k+1}\geq m_{k}$ and $n_{k+1} > n_{k}$ for each positive integer $k$,
\item $m_k\geq n_k$ for each positive integer $k$,
\end{enumerate}
and there is a sequence  $(H_\ell)$ of upper semicontinuous functions $H_\ell:X_{m_\ell}\multimap Y_{n_\ell}$ such that 
\begin{enumerate}
\item[(a)] the graph $\Gamma(H_k)$ is surjective for each positive integer $k$,
\item[(b)] for each $j\in\{1,2,3,\ldots ,n_k\}$, 
$$
g_{j,n_r}(H_r(x))\subseteq g_{j,n_k}(H_k(f_{m_k,m_r}(x)))
$$
holds for each positive integer $r > k$  and  for each $x\in X_{m_r}$ (see Figure \ref{diagram14}),

\begin{figure}[h]
\centering
\begin{tikzpicture}[node distance=1.5cm, auto]
  \node (X1) {$X_{m_k}$};
  \node (X2) [right of=X1] {};
  \node (X3) [right of=X2] {$X_{m_r}$};
  \draw[<-] (X1) to node {$f_{m_k,m_r}$} (X3);
   \node (Z) [below of=X1] {};
    \node (Y3) [below of=Z] {$Y_{n_k}$};
      \node (Y2) [left of=Y3] {};
        \node (Y1) [left of=Y2] {$Y_{j}$};
          \node (Y4) [right of=Y3] {};
            \node (Y5) [right of=Y4] {$Y_{n_r}$};
             \draw[<-] (Y1) to node {$g_{j,n_k}$} (Y3);
              \draw[<-] (Y3) to node {$g_{n_k,n_r}$} (Y5);
               \draw[o-] (Y3) to node {$H_k$} (X1);
                \draw[o-] (Y5) to node {$H_r$} (X3);
\end{tikzpicture}
\caption{The diagram from item $(b)$} \label{diagram14}
\end{figure}
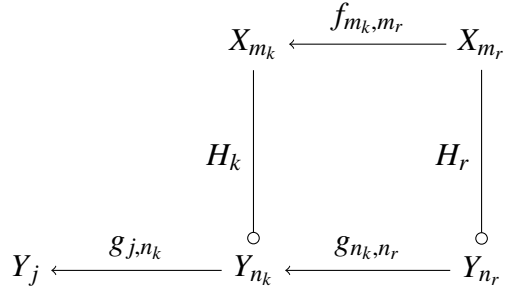

\item[(c)]  for each positive integer $j$, 
$$
\lim_{k\to\infty}\diam (g_{j,n_k}(H_{k}(p_{m_k}(\mathbf x))))=0
$$
for any $\mathbf x\in \varprojlim\{X_\ell,f_\ell\}_{\ell=1}^{\infty}$. 
\end{enumerate}
\end{enumerate}
\end{corollary}
\begin{proof}
The implication from \ref{2ena2} to \ref{1ena1} follows from Theorem \ref{jejhata}. The implication  from \ref{1ena1} to \ref{2ena2} follows from Theorem \ref{Mio:USC2}, which guarantees the existence of the sequences $(m_k)$ and $(n_k)$ of positive integers with the required properties, taking $H_k=T_{n_k,m_k}(g)$ for any continuous surjection $g:\varprojlim\{X_\ell,f_\ell\}_{\ell=1}^{\infty}\rightarrow \varprojlim\{Y_\ell,g_\ell\}_{\ell=1}^{\infty}$.
\end{proof}

The following theorem is our second mapping theorem for inverse limits.

\begin{theorem}\label{mappingTHM2}
Let $\{X_\ell,f_\ell\}_{\ell=1}^{\infty}$ and $\{Y_\ell,g_\ell\}_{\ell=1}^{\infty}$ be inverse sequences of compact metric spaces and surjective continuous bonding functions. The following statements are equivalent.
\begin{enumerate}
    \item\label{1dva1} There is a continuous surjection from $\varprojlim\{X_\ell,f_\ell\}_{\ell=1}^{\infty}$ to $\varprojlim\{Y_\ell,g_\ell\}_{\ell=1}^{\infty}$. 
    \item \label{2dva2} There is a sequence $(H_\ell)$ of upper semicontinuous functions $H_\ell:X_{\ell+1}\multimap Y_{\ell}$ such that 
\begin{enumerate}
\item[(a)]  $g_{k}(H_{k+1}(x))\subseteq H_k(f_{k+1}(x))$
holds for each positive integer $k$  and  for each $x\in X_{k+2}$ (see Figure \ref{diagram15}),

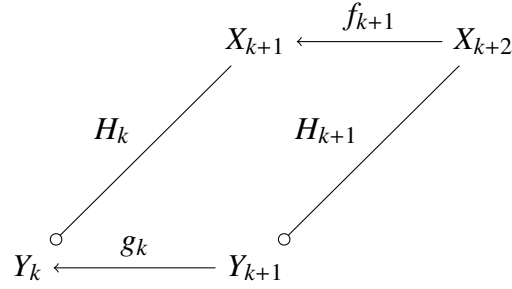
\begin{figure}[h]
\centering
\begin{tikzpicture}[node distance=1.5cm, auto]
  \node (X1) {};
  \node (X2) [right of=X1] {};
  \node (X3) [right of=X2] {$X_{k+1}$};
    \node (X4) [right of=X3] {};
  \node (X5) [right of=X4] {$X_{k+2}$};
  \draw[<-] (X3) to node {$f_{k+1}$} (X5);
   \node (Z) [below of=X1] {};
    \node (Y3) [below of=Z] {$Y_{k}$};
          \node (Y4) [right of=Y3] {};
            \node (Y5) [right of=Y4] {$Y_{k+1}$};
              \draw[<-] (Y3) to node {$g_{k}$} (Y5);
               \draw[o-] (Y3) to node {$H_k$} (X3);
                \draw[o-] (Y5) to node {$H_{k+1}$} (X5);
\end{tikzpicture}
\caption{The diagram from item $(a)$} \label{diagram15}
\end{figure}

\item[(b)]  for each positive integer $j\leq k$, 
$$
\lim_{k\to\infty}\diam (g_{j,k}(H_{k}(p_{k+1}(\mathbf x))))=0 
$$
for any $\mathbf x\in \varprojlim\{X_\ell,f_\ell\}_{\ell=2}^{\infty}$.
\end{enumerate}
\end{enumerate}
\end{theorem}
\begin{proof}
First we prove the implication from \ref{2dva2} to \ref{1dva1}. Let the diagram on Figure \ref{diagram16} be a diagram 
such that all the assumptions from \ref{2dva2} of Theorem \ref{mappingTHM2} are satisfied.

\begin{figure}[h]
\centering
\noindent \begin{tikzpicture}[node distance=1.5cm, auto]
  \node (X1) {};
  \node (X2) [right of=X1] {};
  \node (X3) [right of=X2] {$X_{2}$};
    \node (X4) [right of=X3] {};
  \node (X5) [right of=X4] {$X_{3}$};
   \node (X6) [right of=X5] {};
  \node (X7) [right of=X6] {$X_{4}$};
   \node (X8) [right of=X7] {$\ldots$};
    \node (X9) [right of=X8] {$\varprojlim\{X_\ell,f_\ell\}_{\ell=2}^{\infty}$};
  \draw[<-] (X3) to node {$f_{2}$} (X5);
   \draw[<-] (X5) to node {$f_{3}$} (X7);
    \draw[<-] (X7) to node {} (X8);
   \node (Z1) [below of=X1] {};
   \node (Z2) [right of=Z1] {};
   \node (Z3) [right of=Z2] {};
   \node (Z4) [right of=Z3] {};
   \node (Z5) [right of=Z4] {};
   \node (Z6) [right of=Z5] {};
   \node (Z7) [right of=Z6] {};
   \node (Z8) [right of=Z7] {};
   \node (Z9) [right of=Z8] {};
    \node (Y3) [below of=Z] {$Y_{1}$};
          \node (Y4) [right of=Y3] {};
            \node (Y5) [right of=Y4] {$Y_{2}$};
            \node (Y6) [right of=Y5] {};
            \node (Y7) [right of=Y6] {$Y_{3}$};
            \node (Y8) [right of=Y7] {};
            \node (Y9) [right of=Y8] {$Y_{4}$};
                 \node (Y10) [right of=Y9] {$\ldots$};
     \node (Y11) [right of=Y10] {$\varprojlim\{Y_\ell,g_\ell\}_{\ell=1}^{\infty}$};
              \draw[<-] (Y3) to node {$g_{1}$} (Y5);
               \draw[<-] (Y5) to node {$g_{2}$} (Y7);
                \draw[<-] (Y7) to node {$g_{3}$} (Y9);
               \draw[o-] (Y3) to node {$H_1$} (X3);
                \draw[o-] (Y5) to node {$H_{2}$} (X5);
                 \draw[o-] (Y7) to node {$H_{3}$} (X7);
                 \draw[o-] (Y9) to node {} (Z8);
\end{tikzpicture}
\caption{A diagram for which all the assumptions from \ref{2dva2} of Theorem \ref{mappingTHM2} are satisfied} \label{diagram16}
\end{figure}
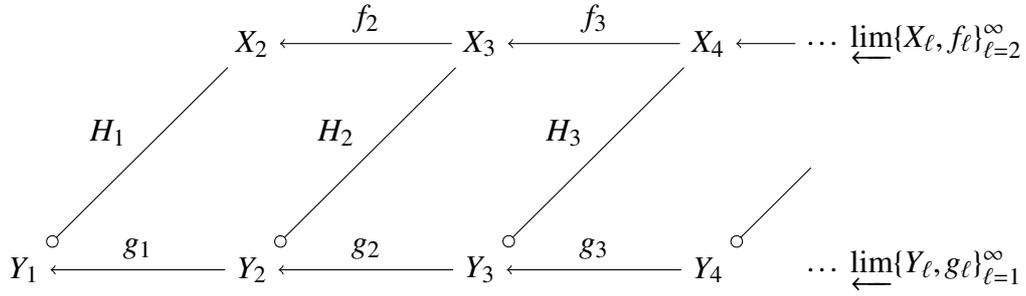
Then, by Theorem \ref{jejhata}, there is a continuous surjection from $\varprojlim\{X_\ell,f_\ell\}_{\ell=1}^{\infty}$ to $\varprojlim\{Y_\ell,g_\ell\}_{\ell=1}^{\infty}$.

Now we prove the implication from \ref{1dva1} to \ref{2dva2}. Let $g:\varprojlim\{X_\ell,f_\ell\}_{\ell=1}^{\infty}\rightarrow \varprojlim\{Y_\ell,g_\ell\}_{\ell=1}^{\infty}$ be a continuous surjection and let 
$m_k=k+1$ and $n_k=k$, and let $H_k=T_{k,k+1}(g)=q_k\circ g\circ p_{k+1}^{-1}$ for each positive integer $k$. Let $k$ be a positive integer and let $x\in X_{k+1}$. Since $p_{k+2}^{-1}(x)\subseteq p_{k+1}^{-1}(f_{k+1}(x))$, it follows from
$$
g_{k}(H_{k+1}(x))=g_k(q_{k+1}(g(p_{k+2}^{-1}(x))))=q_{k}(g(p_{k+2}^{-1}(x)))
$$
and
$$
H_k(f_{k+1}(x))=q_k(g(p_{k+1}^{-1}(f_{k+1}(x))))
$$
that 
$$
g_{k}(H_{k+1}(x))\subseteq H_k(f_{k+1}(x)).
$$
Let $j$ be a positive integer and let $\mathbf x\in \varprojlim\{X_\ell,f_\ell\}_{\ell=2}^{\infty}$ be any point. We show that 
$$
\lim_{k\to\infty}\diam (g_{j,k}(H_{k}(p_{k+1}(\mathbf x))))=0. 
$$
Let $\varepsilon>0$ and let $\delta >0$ such that for each $A\subseteq \varprojlim\{X_\ell,f_\ell\}_{\ell=2}^{\infty}$,
$$
\diam(A)<\delta  \Longrightarrow  \diam(q_j(g(A)))<\varepsilon
 $$
holds. Such a $\delta$ does exist since $q_j\circ g$ is uniformly continuous. Let $k$ be a positive integer such that $\frac{1}{2^k}<\delta$. Then,
$$
\diam (g_{j,k}(H_{k}(p_{k+1}(\mathbf x))))=\diam(q_j(g(p_{k+1}^{-1}(p_{k+1}(\mathbf x)))))<\varepsilon,
$$
since $\diam(p_{k+1}^{-1}(p_{k+1}(\mathbf x)))<\frac{1}{2^k}<\delta$.
This completes the proof. 
\end{proof}
In the rest of this section, we study the fixed point property of inverse limits. 
\begin{theorem}\label{main:1}
Let $\{X_\ell,f_\ell\}_{\ell=1}^{\infty}$ be an inverse sequence of compact metric spaces and surjective continuous bonding functions, and let $g:\varprojlim\{X_\ell,f_\ell\}_{\ell=1}^{\infty}\rightarrow \varprojlim\{X_\ell,f_\ell\}_{\ell=1}^{\infty}$ be a continuous function. 
The following statements are equivalent. 
\begin{enumerate}
\item\label{encka} There is a point $\mathbf x\in \varprojlim\{X_\ell,f_\ell\}_{\ell=1}^{\infty}$ such that $g(\mathbf x)=\mathbf x$.
\item \label{dvencka} 
There are sequences $(n_\ell)$ and $(m_\ell)$ of positive integers 
such that 
\begin{enumerate}
\item[(i)] $m_{k+1}\geq m_{k}$ and $n_{k+1} > n_{k}$ for each positive integer $k$,
\item[(ii)] $m_k\geq n_k$ for each positive integer $k$,
\end{enumerate}
and there is a sequence $(H_\ell)$ of upper semicontinuous functions $H_\ell:X_{m_\ell}\multimap X_{n_\ell}$ such that 
\begin{enumerate}
\item[(a)] \label{uno13} for each positive integer $k$ and for each $j\in\{1,2,3,\ldots ,n_k\}$,
$$
f_{j,n_r}(H_r(x))\subseteq f_{j,n_k}(H_k(f_{m_k,m_r}(x)))
$$
holds for each positive integer $r > k$  and  for each $x\in X_{m_r}$ (see Figure \ref{diagram17}),

\begin{figure}[h]
\centering
\begin{tikzpicture}[node distance=1.5cm, auto]
  \node (X1) {$X_{m_k}$};
  \node (X2) [right of=X1] {};
  \node (X3) [right of=X2] {$X_{m_r}$};
  \draw[<-] (X1) to node {$f_{m_k,m_r}$} (X3);
   \node (Z) [below of=X1] {};
    \node (Y3) [below of=Z] {$X_{n_k}$};
      \node (Y2) [left of=Y3] {};
        \node (Y1) [left of=Y2] {$X_{j}$};
          \node (Y4) [right of=Y3] {};
            \node (Y5) [right of=Y4] {$X_{n_r}$};
             \draw[<-] (Y1) to node {$f_{j,n_k}$} (Y3);
              \draw[<-] (Y3) to node {$f_{n_k,n_r}$} (Y5);
               \draw[o-] (Y3) to node {$H_k$} (X1);
                \draw[o-] (Y5) to node {$H_r$} (X3);
\end{tikzpicture}
\caption{The diagram from item $(a)$} \label{diagram17}
\end{figure}
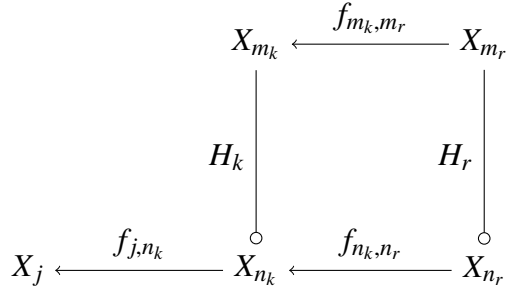

\item[(b)] \label{due23a} for each positive integer $k$ and for each $j\in\{1,2,3,\ldots ,n_k\}$,
$$
p_{j}(g(\mathbf x))\in f_{j,n_k}(H_k(p_{m_k}(\mathbf x)))
$$
holds for each $\mathbf x\in \varprojlim\{X_\ell,f_\ell\}_{\ell=1}^{\infty}$ (see Figure \ref{diagram18}), 

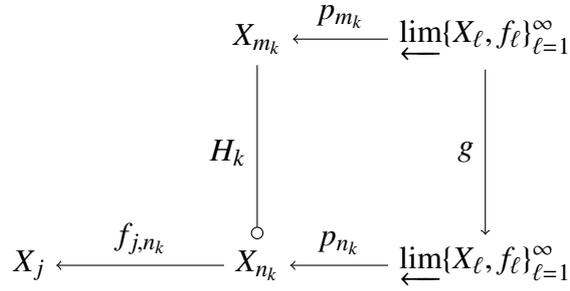
\begin{figure}[h]
\centering
\begin{tikzpicture}[node distance=1.5cm, auto]
  \node (X1) {$X_{m_k}$};
  \node (X2) [right of=X1] {};
  \node (X3) [right of=X2] {$\varprojlim\{X_\ell,f_\ell\}_{\ell=1}^{\infty}$};
 \draw[<-] (X1) to node {$p_{m_k}$} (X3);
   \node (Z) [below of=X1] {};
    \node (Y3) [below of=Z] {$X_{n_k}$};
      \node (Y2) [left of=Y3] {};
        \node (Y1) [left of=Y2] {$X_{j}$};
          \node (Y4) [right of=Y3] {};
           \node (Y5) [right of=Y4] {$\varprojlim\{X_\ell,f_\ell\}_{\ell=1}^{\infty}$};
             \draw[<-] (Y1) to node {$f_{j,n_k}$} (Y3);
              \draw[<-] (Y3) to node {$p_{n_k}$} (Y5);
               \draw[o-] (Y3) to node {$H_k$} (X1);
                \draw[<-] (Y5) to node {$g$} (X3);
\end{tikzpicture}
\caption{The diagram from item $(b)$} \label{diagram18}
\end{figure}

\item[(c)] \label{tre33c} for each positive integer $j$,
$$
\lim_{k\to \infty}\diam (f_{j,n_k}(H_{k}(p_{m_k}(\mathbf x))))=0 
$$
for any $\mathbf x\in \varprojlim\{X_\ell,f_\ell\}_{\ell=1}^{\infty}$,
\item[(d)] \label{tre33d} for each positive integer $k$ and for each $j\in\{1,2,3,\ldots ,n_k\}$, there is a point $x\in X_{m_k}$ such that 
$$
f_{j,m_k}(x)\in f_{j,n_k}(H_k(x))).
$$
\end{enumerate}
\end{enumerate}
\end{theorem}
\begin{proof}
We first prove that \ref{dvencka} follows from \ref{encka}. Let $\mathbf x_0\in \varprojlim\{X_\ell,f_\ell\}_{\ell=1}^{\infty}$ such that $g(\mathbf x_0)=\mathbf x_0$.
It follows from Theorem \ref{Mio:USC2} that 
there are  sequences $(n_\ell)$ and $(m_\ell)$ of positive integers such that 
\begin{enumerate}
\item $m_{k+1}\geq m_{k}$ and $n_{k+1} > n_{k}$ for each positive integer $k$,
\item $m_k\geq n_k$ for each positive integer $k$,
\item for each positive integer $k$ and for each $j\in\{1,2,3,\ldots ,n_k\}$, 
$$
f_{j,n_r}(T_{n_r,m_r}(g)(x))\subseteq f_{j,n_k}(T_{n_k,m_k}(g)(f_{m_k,m_r}(x)))
$$
holds for each positive integer $r > k$  and  for each $x\in X_{m_r}$,
\item for each positive integer $k$ and for each $j\in\{1,2,3,\ldots ,n_k\}$,  
$$
p_{j}(g(\mathbf x))\in f_{j,n_k}(T_{n_k,m_k}(g)(p_{m_k}(\mathbf x)))
$$
holds for each $\mathbf x\in \varprojlim\{X_\ell,f_\ell\}_{\ell=1}^{\infty}$,
\item for each positive integer $j$,
$$
\lim_{k\to\infty}\diam (f_{j,n_k}(T_{n_k,m_k}(g)(p_{m_k}(\mathbf x))))=0
$$
for any $\mathbf x\in \varprojlim\{X_\ell,f_\ell\}_{\ell=1}^{\infty}$.
\end{enumerate}
Choose and fix such sequences $(n_\ell)$ and $(m_\ell)$. Set $H_k=T_{n_k,m_k}(g)$ for each positive integer $k$. All that is left to prove is that for each positive integer $k$ and for each $j\in\{1,2,3,\ldots ,n_k\}$, there is a point $x\in X_{m_k}$ such that 
$$
f_{j,m_k}(x)\in f_{j,n_k}(H_k(x)).
$$
Let $k$ be a positive integer and let $j\in\{1,2,3,\ldots ,n_k\}$. Let $x=p_{m_k}(\mathbf x_0)$. Then 
$$
f_{j,m_k}(x)=f_{j,m_k}(p_{m_k}(\mathbf x_0))=p_j(\mathbf x_0)
$$
and
$$
f_{j,n_k}(H_k(x)))=f_{j,n_k}(T_{n_k,m_k}(g)(x)))=f_{j,n_k}(p_{n_k}(g(p_{m_k}^{-1}(p_{m_k}(\mathbf x_0))))).
$$
Since $\mathbf x_0\in p_{m_k}^{-1}(p_{m_k}(\mathbf x_0))$, it follows that
$$
f_{j,n_k}(p_{n_k}(g(\mathbf x_0)))\in f_{j,n_k}(p_{n_k}(g(p_{m_k}^{-1}(p_{m_k}(\mathbf x_0))))).
$$
It follows from  $f_{j,n_k}(p_{n_k}(g(\mathbf x_0)))=f_{j,n_k}(p_{n_k}(\mathbf x_0))=p_j(\mathbf x_0)$ that
$$
f_{j,m_k}(x)\in f_{j,n_k}(H_k(x)).
$$
Therefore, \ref{dvencka} follows. 

Next, we prove that \ref{encka} follows from \ref{dvencka}. 
Let the diagram on Figure \ref{diagram19} be a diagram such that all the assumptions of \ref{dvencka} are satisfied.

\begin{figure}[h]
\centering
\noindent \begin{tikzpicture}[node distance=1.5cm, auto]
  \node (X1) {$X_{m_1}$};
  \node (X2) [right of=X1] {};
   \node (X3) [right of=X2] {$X_{m_2}$};
    \node (X4) [right of=X3] {};
      \node (X5) [right of=X4] {$X_{m_3}$};
    \node (X6) [right of=X5] {};
        \node (X7) [right of=X6] {$X_{m_4}$};
    \node (X8) [right of=X7] {$\ldots$};
    \node (X9) [right of=X8] {$\varprojlim\{X_\ell,f_\ell\}_{\ell=1}^{\infty}$};
   \node (Z1) [below of=X1] {};
  \node (Z2) [right of=Z1] {};
    \node (Z3) [right of=Z2] {};
    \node (Z4) [right of=Z3] {};
      \node (Z5) [right of=Z4] {};
    \node (Z6) [right of=Z5] {};
        \node (Z7) [right of=Z6] {};
    \node (Z8) [right of=Z7] {};
    \node (Z9) [right of=Z8] {};
  \draw[<-] (X1) to node {$f_{m_1,m_2}$} (X3);
   \draw[<-] (X3) to node {$f_{m_2,m_3}$} (X5);
      \draw[<-] (X5) to node {$f_{m_3,m_4}$} (X7);
          \draw[<-] (X7) to node {$f_{m_4}$} (X8);
\node (Y1) [below of=Z1] {$X_{n_1}$};
  \node (Y2) [right of=Y1] {};
    \node (Y3) [right of=Y2] {$X_{n_2}$};
      \node (a) [right of=Y3] {};
       \node (b) [right of=a] {$X_{n_3}$};
     \node (Y4) [right of=b] {};
      \node (Y5) [right of=Y4] {$X_{n_4}$};
         \node (Y6) [right of=Y5] {$\ldots$};
      \node (Y7) [below of=Z9] {$\varprojlim\{X_\ell,f_\ell\}_{\ell=1}^{\infty}$};
  \draw[<-] (Y1) to node {$f_{n_1,n_2}$} (Y3);
   \draw[<-] (Y3) to node {$f_{n_2,n_3}$} (b);
     \draw[<-] (b) to node {$f_{n_3,n_4}$} (Y5);
   \draw[<-] (Y5) to node {$f_{n_4}$} (Y6);
      \draw[o-] (Y1) to node {$H_{1}$} (X1);
            \draw[o-] (Y3) to node {$H_{2}$} (X3);
                  \draw[o-] (b) to node {$H_{3}$} (X5);
                        \draw[o-] (Y5) to node {$H_{4}$} (X7);
      \draw[<-] (Y7) to node {$g$} (X9);
\end{tikzpicture}
\caption{A diagram such that all the assumptions of \ref{dvencka} are satisfied} \label{diagram19}
\end{figure}
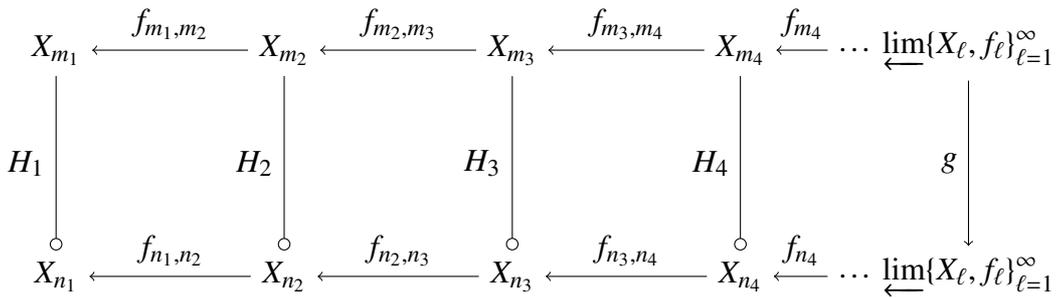
For each positive integer $k$, let $x_{m_k}\in X_{m_k}$ be such a point that 
$$
f_{j,m_k}(x_{m_k})\in f_{j,n_k}(H_k(x_{m_k}))
$$
for each $j\in \{1,2,3,\ldots ,n_k\}$. Also, for each positive integer $k$, let
$$
\mathbf{x}^{m_k}=(x_1^{m_k},x_2^{m_k},x_3^{m_k},\ldots )\in \varprojlim\{X_\ell,f_\ell\}_{\ell=1}^{\infty}
$$
such that $p_{m_k}(\mathbf{x}^{m_k})=x_{m_k}$. Let $\mathbf x^0\in \varprojlim\{X_\ell,f_\ell\}_{\ell=1}^{\infty}$ be an accumulation point of the sequence $(\mathbf{x}^{m_k})$. Without any loss of generality, we assume that 
$$
\lim_{k\to \infty} \mathbf x^{m_{k}}=\mathbf x^0.
$$
We show that 
$$
g(\mathbf x^0)=\mathbf x^0
$$
by showing that for each positive integer $j$,
$$
p_{j}(g(\mathbf x^0))=p_{j}(\mathbf x^0).
$$
Let $j$ be a positive integer. Recall that by Observation \ref{tupki}, 
$$
\bigcap_{k=1}^{\infty}f_{j,n_{k}}(H_{k}(p_{m_{k}}(\mathbf x^0)))=\{p_{j}(g(\mathbf x^0))\}.
$$
Note also that for each positive integer $k$,
$$
f_{j,n_{k+1}}(H_{k+1}(p_{m_{k+1}}(\mathbf x^{0})))\subseteq f_{j,n_{k}}(H_{k}(f_{m_{k},m_{k+1}}(p_{m_{k+1}}(\mathbf x^{0})))),
$$
and therefore,
$$
p_{j}(g(\mathbf x^0)) \in \bigcap_{k=1}^{\infty}f_{j,n_{k}}(H_{k}(f_{m_{k},m_{k+1}}(p_{m_{k+1}}(\mathbf x^{0})))).
$$
Since 
$$
\lim_{k\to\infty}\diam (f_{j,n_{k}}(H_{k}(f_{m_{k},m_{k+1}}(p_{m_{k+1}}(\mathbf x^{0})))))=0,
$$
it follows that 
$$
\diam (\bigcap_{k=1}^{\infty}f_{j,n_{k}}(H_{k}(f_{m_{k},m_{k+1}}(p_{m_{k+1}}(\mathbf x^{0})))))=0. 
$$
Therefore,
$$
\bigcap_{k=1}^{\infty}f_{j,n_{k}}(H_{k}(f_{m_{k},m_{k+1}}(p_{m_{k+1}}(\mathbf x^{0}))))=\{p_{j}(g(\mathbf x^0)) \}.
$$
Since $f_{m_{k},m_{k+1}}(p_{m_{k+1}}(\mathbf x^{0}))=p_{m_{k}}(\mathbf x^{0})$, it follows that 
$$
\bigcap_{k=1}^{\infty}f_{j,n_{k}}(H_{k}(f_{m_{k},m_{k+1}}(p_{m_{k+1}}(\mathbf x^{0}))))=\bigcap_{k=1}^{\infty}f_{j,n_{k}}(H_{k}(p_{m_{k}}(\mathbf x^{0}))).
$$
On the other hand, also taking into account that each $H_\ell$ is upper semicontinuous, we get that for each positive integer $k$ such that $j\leq n_k$,
$$
p_{j}(\mathbf x^0)=f_{j,m_{k}}(p_{m_{k}}(\mathbf x^{0}))=
f_{j,m_{k}}(p_{m_{k}}(\lim_{\ell\to \infty}\mathbf x^{m_\ell}))=
$$
$$
\lim_{\ell\to \infty}f_{j,m_{k}}(p_{m_{k}}(\mathbf x^{m_\ell}))=
\lim_{\ell\to \infty}f_{j,m_{k}}(f_{m_{k},m_{\ell}}(p_{m_\ell}(\mathbf x^{m_\ell}))=
$$
$$
\lim_{\ell\to \infty}f_{j,m_{\ell}}(x_{m_\ell})\in \lim_{\ell\to \infty}f_{j,n_{\ell}}(H_{{\ell}}(x_{m_\ell}))
$$
with respect to the Hausdorff metric. Note that $\displaystyle \lim_{\ell\to \infty}f_{j,n_{\ell}}(H_{{\ell}}(x_{m_\ell}))$ does exist and that 
$$
\lim_{\ell\to \infty}f_{j,n_{\ell}}(H_{{\ell}}(x_{m_\ell}))=\bigcap_{\ell=1}^{\infty}f_{j,n_{\ell}}(H_{{\ell}}(x_{m_\ell})).
$$
Note also that
$$
\lim_{\ell\to \infty}f_{j,n_{\ell}}(H_{{\ell}}(x_{m_\ell}))=\lim_{\ell\to \infty}f_{j,n_{\ell}}(H_{{\ell}}(p_{m_\ell}(\mathbf x^{m_\ell})))\subseteq 
$$
$$
\lim_{\ell\to \infty}f_{j,n_{k}}(H_{{k}}(f_{m_{k},m_{\ell}}(p_{m_\ell}(\mathbf x^{m_\ell}))))=\lim_{\ell\to \infty}f_{j,n_{k}}(H_{{k}}(p_{m_{k}}(\mathbf x^{m_\ell})))=
$$
$$
f_{j,n_{k}}(H_{{k}}(p_{m_{k}}(
\lim_{\ell\to \infty}\mathbf x^{m_\ell})))=f_{j,n_{k}}(H_{k}(p_{m_{k}}(\mathbf x^{0}))).
$$
 We have proved that for each positive integer $k$ such that $j\leq n_k$,
 $$
 p_{j}(\mathbf x^0)\in f_{j,n_{k}}(H_{k}(p_{m_{k}}(\mathbf x^{0})))
 $$
 and therefore, 
$$
p_{j}(\mathbf x^0)\in \bigcap_{k=1}^{\infty}f_{j,n_{k}}(H_{k}(p_{m_{k}}(\mathbf x^{0}))).
$$
It follows that $p_{j}(\mathbf x^0)=p_{j}(g(\mathbf x^0))$. 
\end{proof}
\begin{corollary}\label{main:2}
Let $\{X_\ell,f_\ell\}_{\ell=1}^{\infty}$ be an inverse sequence of compact metric spaces and surjective continuous bonding functions, and let $g:\varprojlim\{X_\ell,f_\ell\}_{\ell=1}^{\infty}\rightarrow \varprojlim\{X_\ell,f_\ell\}_{\ell=1}^{\infty}$ be a continuous function.

The following statements are equivalent. 
\begin{enumerate}
\item For each  point $\mathbf x\in \varprojlim\{X_\ell,f_\ell\}_{\ell=1}^{\infty}$,  
$$
g(\mathbf x)\neq \mathbf x.
$$
\item For all  sequences $(n_\ell)$ and $(m_\ell)$ of positive integers 
such that 
\begin{enumerate}
\item[(i)] $m_{k+1}\geq m_{k}$ and $n_{k+1} > n_{k}$ for each positive integer $k$,
\item[(ii)] $m_k\geq n_k$ for each positive integer $k$,
\end{enumerate}
and for each sequence $(H_\ell)$ of upper semicontinuous functions $H_\ell:X_{m_\ell}\multimap X_{n_\ell}$ such that 
\begin{enumerate}
\item[(a)] \label{uno13b} for each positive integer $k$ and for each $j\in\{1,2,3,\ldots ,n_k\}$, 
$$
f_{j,n_r}(H_r(x))\subseteq f_{j,n_k}(H_k(f_{m_k,m_r}(x)))
$$
holds for each positive integer $r > k$  and  for each $x\in X_{m_r}$ (see Figure \ref{diagram20}),

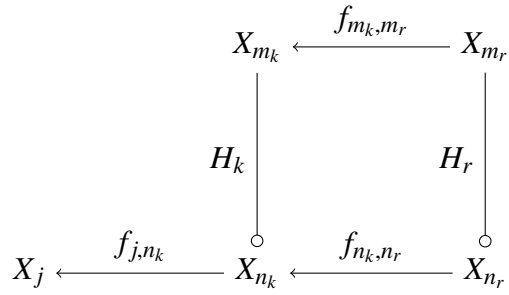
\begin{figure}[h]
\centering
\begin{tikzpicture}[node distance=1.5cm, auto]
  \node (X1) {$X_{m_k}$};
  \node (X2) [right of=X1] {};
  \node (X3) [right of=X2] {$X_{m_r}$};
  \draw[<-] (X1) to node {$f_{m_k,m_r}$} (X3);
   \node (Z) [below of=X1] {};
    \node (Y3) [below of=Z] {$X_{n_k}$};
      \node (Y2) [left of=Y3] {};
        \node (Y1) [left of=Y2] {$X_{j}$};
          \node (Y4) [right of=Y3] {};
            \node (Y5) [right of=Y4] {$X_{n_r}$};
             \draw[<-] (Y1) to node {$f_{j,n_k}$} (Y3);
              \draw[<-] (Y3) to node {$f_{n_k,n_r}$} (Y5);
               \draw[o-] (Y3) to node {$H_k$} (X1);
                \draw[o-] (Y5) to node {$H_r$} (X3);
\end{tikzpicture}
\caption{The diagram from item $(a)$} \label{diagram20}
\end{figure}

\item[(b)] \label{due23} for each positive integer $k$ and for each $j\in\{1,2,3,\ldots ,n_k\}$, 
$$
f_{j,n_k}(p_{n_k}(g(\mathbf x)))\in f_{j,n_k}(H_k(p_{m_k}(\mathbf x)))
$$
holds for each $\mathbf x\in \varprojlim\{X_\ell,f_\ell\}_{\ell=1}^{\infty}$ (see Figure \ref{diagram21}), 

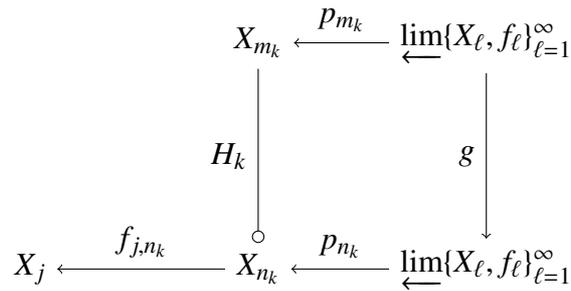
\begin{figure}[h]
\centering
\begin{tikzpicture}[node distance=1.5cm, auto]
  \node (X1) {$X_{m_k}$};
  \node (X2) [right of=X1] {};
  \node (X3) [right of=X2] {$\varprojlim\{X_\ell,f_\ell\}_{\ell=1}^{\infty}$};
 \draw[<-] (X1) to node {$p_{m_k}$} (X3);
   \node (Z) [below of=X1] {};
    \node (Y3) [below of=Z] {$X_{n_k}$};
      \node (Y2) [left of=Y3] {};
        \node (Y1) [left of=Y2] {$X_{j}$};
          \node (Y4) [right of=Y3] {};
           \node (Y5) [right of=Y4] {$\varprojlim\{X_\ell,f_\ell\}_{\ell=1}^{\infty}$};
             \draw[<-] (Y1) to node {$f_{j,n_k}$} (Y3);
              \draw[<-] (Y3) to node {$p_{n_k}$} (Y5);
               \draw[o-] (Y3) to node {$H_k$} (X1);
                \draw[<-] (Y5) to node {$g$} (X3);
\end{tikzpicture}
\caption{The diagram from item $(b)$} \label{diagram21}
\end{figure}

\item[(c)] \label{tre33} for each positive integer $j$,
$$
\lim_{k\to \infty}\diam (f_{j,n_k}(H_{k}(p_{m_k}(\mathbf x))))=0 
$$
for any $\mathbf x\in \varprojlim\{X_\ell,f_\ell\}_{\ell=1}^{\infty}$,
\end{enumerate}
there are positive integers $k$ and $j\in \{1,2,3,\ldots ,n_k\}$ such that for each point $x\in X_{m_k}$,   
$$
f_{j,m_k}(x)\not \in f_{j,n_k}(H_k(x)).
$$
\end{enumerate}
\end{corollary}
\begin{proof}
The corollary follows from Theorem \ref{main:1}. 
\end{proof}
We conclude the section by stating and proving the following theorem.
\begin{theorem}\label{main:3}
Let $\{X_\ell,f_\ell\}_{\ell=1}^{\infty}$ be an inverse sequence of compact metric spaces and surjective continuous bonding functions.
The following statements are equivalent. 
\begin{enumerate}
\item\label{muka} The inverse limit $\varprojlim\{X_\ell,f_\ell\}_{\ell=1}^{\infty}$ does not have the fixed point property.
\item\label{maka} There is a sequence $(H_\ell)$ of upper semicontinuous functions $H_\ell:X_{\ell+1}\multimap X_{\ell}$ such that 
\begin{enumerate}
\item[(a)] \label{jed} $f_{k}(H_{k+1}(x))\subseteq H_k(f_{k+1}(x))$
holds for each positive integer $k$  and  for each $x\in X_{k+1}$ (see Figure \ref{diagram22}),

\begin{figure}[h]
\centering
\begin{tikzpicture}[node distance=1.5cm, auto]
  \node (X1) {};
  \node (X2) [right of=X1] {};
  \node (X3) [right of=X2] {$X_{k+1}$};
    \node (X4) [right of=X3] {};
  \node (X5) [right of=X4] {$X_{k+2}$};
  \draw[<-] (X3) to node {$f_{k+1}$} (X5);
   \node (Z) [below of=X1] {};
    \node (Y3) [below of=Z] {$X_{k}$};
          \node (Y4) [right of=Y3] {};
            \node (Y5) [right of=Y4] {$X_{k+1}$};
              \draw[<-] (Y3) to node {$f_{k}$} (Y5);
               \draw[o-] (Y3) to node {$H_k$} (X3);
                \draw[o-] (Y5) to node {$H_{k+1}$} (X5);
\end{tikzpicture}
\caption{The diagram from item $(a)$} \label{diagram22}
\end{figure}
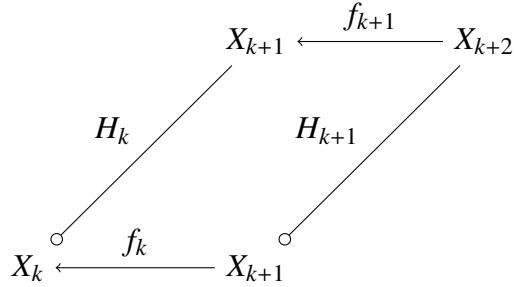

\item[(b)] \label{dv} for each positive integer $j$, 
$$
\lim_{k\to\infty}\diam (f_{j,k}(H_{k}(p_{k+1}(\mathbf x))))=0 
$$
for any $\mathbf x\in \varprojlim\{X_\ell,f_\ell\}_{\ell=2}^{\infty}$,
\item[(c)] \label{tr} there is a positive integer $k$, such that for each point $x\in X_{k+1}$,
$$
f_{k}(x)\not \in H_k(x).
$$
\end{enumerate}
\end{enumerate}
\end{theorem}
\begin{proof}
First we prove the implication from \ref{maka} to \ref{muka}. Let the diagram on Figure \ref{diagram23} be a diagram such that all the assumptions of \ref{maka} are satisfied.

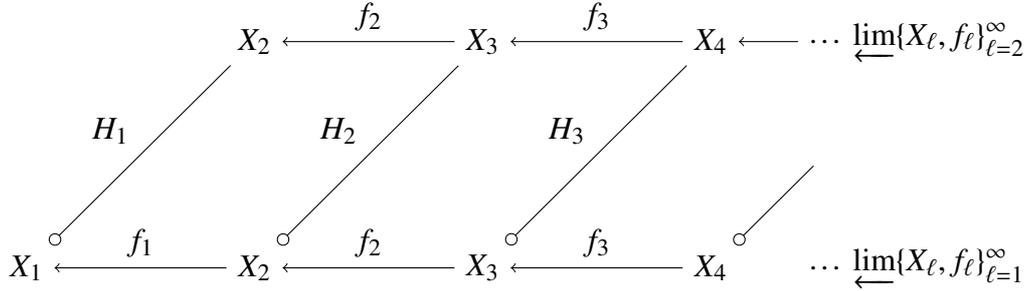
\begin{figure}[h]
\centering
\noindent \begin{tikzpicture}[node distance=1.5cm, auto]
  \node (X1) {};
  \node (X2) [right of=X1] {};
  \node (X3) [right of=X2] {$X_{2}$};
    \node (X4) [right of=X3] {};
  \node (X5) [right of=X4] {$X_{3}$};
   \node (X6) [right of=X5] {};
  \node (X7) [right of=X6] {$X_{4}$};
   \node (X8) [right of=X7] {$\ldots$};
    \node (X9) [right of=X8] {$\varprojlim\{X_\ell,f_\ell\}_{\ell=2}^{\infty}$};
  \draw[<-] (X3) to node {$f_{2}$} (X5);
   \draw[<-] (X5) to node {$f_{3}$} (X7);
    \draw[<-] (X7) to node {} (X8);
   \node (Z1) [below of=X1] {};
   \node (Z2) [right of=Z1] {};
   \node (Z3) [right of=Z2] {};
   \node (Z4) [right of=Z3] {};
   \node (Z5) [right of=Z4] {};
   \node (Z6) [right of=Z5] {};
   \node (Z7) [right of=Z6] {};
   \node (Z8) [right of=Z7] {};
   \node (Z9) [right of=Z8] {};
    \node (Y3) [below of=Z] {$X_{1}$};
          \node (Y4) [right of=Y3] {};
            \node (Y5) [right of=Y4] {$X_{2}$};
            \node (Y6) [right of=Y5] {};
            \node (Y7) [right of=Y6] {$X_{3}$};
            \node (Y8) [right of=Y7] {};
            \node (Y9) [right of=Y8] {$X_{4}$};
                 \node (Y10) [right of=Y9] {$\ldots$};
     \node (Y11) [right of=Y10] {$\varprojlim\{X_\ell,f_\ell\}_{\ell=1}^{\infty}$};
              \draw[<-] (Y3) to node {$f_{1}$} (Y5);
               \draw[<-] (Y5) to node {$f_{2}$} (Y7);
                \draw[<-] (Y7) to node {$f_{3}$} (Y9);
               \draw[o-] (Y3) to node {$H_1$} (X3);
                \draw[o-] (Y5) to node {$H_{2}$} (X5);
                 \draw[o-] (Y7) to node {$H_{3}$} (X7);
                 \draw[o-] (Y9) to node {} (Z8);
\end{tikzpicture}
\caption{A diagram such that all the assumptions of \ref{maka} are satisfied} \label{diagram23}
\end{figure}
Let 
$$
g:\varprojlim\{X_\ell,f_\ell\}_{\ell=1}^{\infty}\rightarrow \varprojlim\{X_\ell,f_\ell\}_{\ell=1}^{\infty}
$$
be a continuous function such that for any positive integer $k$,
$$
p_k(g(\mathbf x))\in H_k(p_{k+1}(\mathbf x))
$$
for any $\mathbf x\in \varprojlim\{X_\ell,f_\ell\}_{\ell=1}^{\infty}$. Such a function does exist by Theorem \ref{jejhata}.  We show that 
$$
g(\mathbf x)\neq \mathbf x
$$
for any
$\mathbf x\in \varprojlim\{X_\ell,f_\ell\}_{\ell=1}^{\infty}$.
Suppose that there is $\mathbf x_0\in \varprojlim\{X_\ell,f_\ell\}_{\ell=1}^{\infty}$ such that $g(\mathbf x_0)=\mathbf x_0$. Fix such a point $\mathbf x_0$. Let $k$ be a positive integer such that $f_k(x)\not \in H_k(x)$ for each $x\in X_{k+1}$. Let $x=p_{k+1}(\mathbf x_0)$. Then 
$$
f_k(x)=f_k(p_{k+1}(\mathbf x_0))=f_k(p_{k+1}(g(\mathbf x_0)))\in f_k(H_{k+1}(p_{k+2}(\mathbf x_0)))\subseteq 
$$
$$
H_k(f_{k+1}(p_{k+2}(\mathbf x_0)))=H_k(p_{k+1}(\mathbf x_0))=H_k(x),
$$
which is a contradiction. Therefore, 
$g(\mathbf x)\neq \mathbf x$ for every $\mathbf x\in \varprojlim\{X_\ell,f_\ell\}_{\ell=1}^{\infty}$ 
and we have proved that the inverse limit $\varprojlim\{X_\ell,f_\ell\}_{\ell=1}^{\infty}$ does not have the fixed point property.

Now we prove the implication from \ref{muka} to \ref{maka}. Let $g:\varprojlim\{X_\ell,f_\ell\}_{\ell=1}^{\infty}\rightarrow \varprojlim\{X_\ell,f_\ell\}_{\ell=1}^{\infty}$ be a continuous function such that $g(\mathbf x)\neq \mathbf x$ for each $\mathbf x\in \varprojlim\{X_\ell,f_\ell\}_{\ell=1}^{\infty}$. 
It follows from Corollary \ref{main:2} that for all  sequences $(n_\ell)$ and $(m_\ell)$ of positive integers 
such that 
\begin{enumerate}
\item[(i)] $m_{k+1}\geq m_{k}$ and $n_{k+1} > n_{k}$ for each positive integer $k$,
\item[(ii)] $m_k\geq n_k$ for each positive integer $k$,
\end{enumerate}
and for each sequence $(H_\ell)$ of upper semicontinuous functions $H_\ell:X_{m_\ell}\multimap X_{n_\ell}$ such that 
\begin{enumerate}
\item[(a)]  for each positive integer $k$ and for each $j\in\{1,2,3,\ldots ,n_k\}$, 
$$
f_{j,n_r}(H_r(x))\subseteq f_{j,n_k}(H_k(f_{m_k,m_r}(x)))
$$
holds for each positive integer $r > k$  and  for each $x\in X_{m_r}$ (see Figure \ref{diagram24}),

\begin{figure}[h]
\centering
\begin{tikzpicture}[node distance=1.5cm, auto]
  \node (X1) {$X_{m_k}$};
  \node (X2) [right of=X1] {};
  \node (X3) [right of=X2] {$X_{m_r}$};
  \draw[<-] (X1) to node {$f_{m_k,m_r}$} (X3);
   \node (Z) [below of=X1] {};
    \node (Y3) [below of=Z] {$X_{n_k}$};
      \node (Y2) [left of=Y3] {};
        \node (Y1) [left of=Y2] {$X_{j}$};
          \node (Y4) [right of=Y3] {};
            \node (Y5) [right of=Y4] {$X_{n_r}$};
             \draw[<-] (Y1) to node {$f_{j,n_k}$} (Y3);
              \draw[<-] (Y3) to node {$f_{n_k,n_r}$} (Y5);
               \draw[o-] (Y3) to node {$H_k$} (X1);
                \draw[o-] (Y5) to node {$H_r$} (X3);
\end{tikzpicture}
\caption{The diagram from item $(a)$} \label{diagram24}
\end{figure}
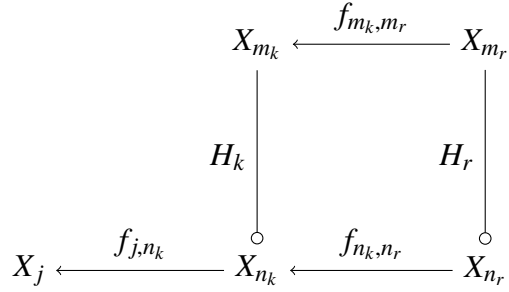

\item[(b)]  for each positive integer $k$ and for each $j\in\{1,2,3,\ldots ,n_k\}$, 
$$
f_{j,n_k}(p_{n_k}(g(\mathbf x)))\in f_{j,n_k}(H_k(p_{m_k}(\mathbf x)))
$$
holds for each $\mathbf x\in \varprojlim\{X_\ell,f_\ell\}_{\ell=1}^{\infty}$ (see Figure \ref{diagram25}), 

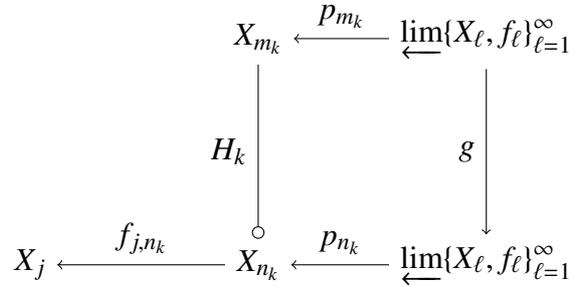
\begin{figure}[h]
\centering
\begin{tikzpicture}[node distance=1.5cm, auto]
  \node (X1) {$X_{m_k}$};
  \node (X2) [right of=X1] {};
  \node (X3) [right of=X2] {$\varprojlim\{X_\ell,f_\ell\}_{\ell=1}^{\infty}$};
 \draw[<-] (X1) to node {$p_{m_k}$} (X3);
   \node (Z) [below of=X1] {};
    \node (Y3) [below of=Z] {$X_{n_k}$};
      \node (Y2) [left of=Y3] {};
        \node (Y1) [left of=Y2] {$X_{j}$};
          \node (Y4) [right of=Y3] {};
           \node (Y5) [right of=Y4] {$\varprojlim\{X_\ell,f_\ell\}_{\ell=1}^{\infty}$};
             \draw[<-] (Y1) to node {$f_{j,n_k}$} (Y3);
              \draw[<-] (Y3) to node {$p_{n_k}$} (Y5);
               \draw[o-] (Y3) to node {$H_k$} (X1);
                \draw[<-] (Y5) to node {$g$} (X3);
\end{tikzpicture}
\caption{The diagram from item $(b)$} \label{diagram25}
\end{figure}

\item[(c)]  for each positive integer $j$,
$$
\lim_{k\to \infty}\diam (f_{j,n_k}(H_{k}(p_{m_k}(\mathbf x))))=0 
$$
for any $\mathbf x\in \varprojlim\{X_\ell,f_\ell\}_{\ell=1}^{\infty}$,
\end{enumerate}
there are positive integers $k$ and $j\in \{1,2,3,\ldots ,n_k\}$ such that for each point $x\in X_{m_k}$,   
$$
f_{j,m_k}(x)\not \in f_{j,n_k}(H_k(x)).
$$
Let $m_k=k+1$ and $n_k=k$, and let $H_k=T_{k,k+1}(g)=p_k\circ g\circ p_{k+1}^{-1}$ for each positive integer $k$. Also, let $k_0$ be a positive integer and $j_0\in \{1,2,3,\ldots,n_{k_0}\}$ such that for each point $x\in X_{m_{k_0}}$,   
$$
f_{j_0,m_{k_0}}(x)\not \in f_{j_0,n_{k_0}}(H_{k_0}(x)).
$$
Therefore, \ref{maka} follows. 
This completes the proof. 
\end{proof}
\section{Applications to inverse limits with set-valued bonding functions}\label{s5}
In the last section, we apply the main results from inverse limits with single-valued bonding functions to inverse limits with set-valued bonding functions.
We begin with the following definition. 
\begin{definition}
Let $\{X_\ell,F_\ell\}_{\ell=1}^{\infty}$ be an inverse sequence of compact metric spaces with upper semicontinuous bonding functions. Then 
$$
\star_{i=m}^n\Gamma(F_i^{-1})=\{(x_m,x_{m+1},\ldots,x_n,x_{n+1})\in \prod_{i=m}^{n+1} X_i \ | \ x_i\in F_i(x_{i+1}) \textup{ for all } i\}
$$ 
for all positive integers $m$ and $n$ such that $m\leq n$. The space $\star_{i=m}^n\Gamma(F_i^{-1})$ is called the {\em Mahavier product of $\Gamma(F_{m}^{-1})$, $\Gamma(F_{m+1}^{-1})$, $\Gamma(F_{m+2}^{-1})$, $\ldots$, $\Gamma(F_{n}^{-1})$\/}.
\end{definition}
In this section, we use the following theorem \cite[Theorem 3.2]{van}.
\begin{theorem}\label{vanizrek2}
Let $\{X_\ell,F_\ell\}_{\ell=1}^{\infty}$ be an inverse sequence of compact metric spaces with upper semicontinuous bonding functions, and let $f_\ell:\star_{i=1}^{\ell+1}\Gamma(F_i^{-1})\rightarrow \star_{i=1}^{\ell}\Gamma(F_i^{-1})$ be defined by 
$$
f_\ell(x_1,x_2,x_3,\ldots,x_{\ell+1},x_{\ell+2})=(x_1,x_2,x_3,\ldots,x_{\ell+1})
$$ 
for each positive integer $\ell$. Then the inverse limits 
$$
\varproj\{X_\ell,F_\ell\}_{\ell=1} ^\infty \textup{ ~~~ and ~~~ } \varprojlim\{\star_{i=1}^{\ell}\Gamma(F_i^{-1}), f_\ell\}_{\ell=1} ^\infty
$$ 
are homeomorphic.
\end{theorem}

We continue by stating and proving the following theorem, which is a mapping theorem for inverse limits with set-valued bonding functions.

\begin{theorem}\label{hopsasa}
Let $\{X_\ell,F_\ell\}_{\ell=1}^{\infty}$ and $\{Y_\ell,G_\ell\}_{\ell=1}^{\infty}$ be inverse sequences of compact metric spaces $X_\ell$ and $Y_\ell$, and upper semicontinuous set-valued functions $F_\ell:X_{\ell+1}\multimap X_\ell$ and $G_\ell:Y_{\ell+1}\multimap Y_\ell$.  

The following statements are equivalent.
\begin{enumerate}
\item There is a continuous surjection from $\varproj\{X_\ell,F_\ell\}_{\ell=1} ^\infty$ to $\varproj\{Y_\ell,G_\ell\}_{\ell=1} ^\infty$.
\item  \label{zwei2} There is a sequence $(H_\ell)$ of upper semicontinuous functions 
$$
H_\ell:\star_{i=1}^{\ell+1}\Gamma(F_i^{-1})\multimap \star_{i=1}^{\ell}\Gamma(G_i^{-1})
$$
such that 
\begin{enumerate}
\item[(a)]  for each positive integer $k$, for each  $(x_1,x_2,x_3,\ldots, x_{k+2},x_{k+3})\in \star_{i=1}^{k+2}\Gamma(F_i^{-1})$ and for each $(z_1,z_2,z_3,\ldots ,z_{k+1},z_{k+2})\in H_{k+1}(x_1,x_2,x_3,\ldots, x_{k+2},x_{k+3})$,
$$
(z_1,z_2,z_3,\ldots ,z_{k+1})\in H_{k}(x_1,x_2,x_3,\ldots, x_{k+2}).
$$
\item[(b)]  for each positive integer $j$ and for each $\mathbf x=(x_1,x_2,x_3,\ldots)\in \varproj\{X_\ell,F_\ell\}_{\ell=1} ^\infty$, 
$$
\lim_{k\to\infty}\diam (\mathcal{H}_{j,k}(\mathbf x))=0, 
$$
where for each positive integer $k>j$, $\mathcal{H}_{j,k}(\mathbf x)$ is the set of all points 
$$
(z_1,z_2,z_3,\ldots ,z_{j+1})\in \star_{i=1}^{j}\Gamma(G_i^{-1})
$$
such that there is a point $(y_1,y_2,y_3,\ldots,y_{k+1})\in H_k(x_1,x_2,x_3,\ldots, x_{k+2})$ such that $(y_1,y_2,y_3,\ldots,y_{j+1})=(z_1,z_2,z_3,\ldots,z_{j+1})$.
\end{enumerate}
\end{enumerate}
\end{theorem}
\begin{proof}
Note that it follows from Theorem \ref{vanizrek2} there is a continuous surjection from $\varproj\{X_\ell,F_\ell\}_{\ell=1} ^\infty$ to $\varproj\{Y_\ell,G_\ell\}_{\ell=1} ^\infty$
if and only if there is a continuous surjection from $\varprojlim\{\star_{i=1}^{\ell}\Gamma(F_i^{-1}),f_\ell\}_{\ell=1} ^\infty$ to $\varprojlim\{\star_{i=1}^{\ell}\Gamma(G_i^{-1}),g_\ell\}_{\ell=1} ^\infty$, 
where all functions $f_\ell:\star_{i=1}^{\ell+1}\Gamma(F_i^{-1})\rightarrow \star_{i=1}^{\ell}\Gamma(F_i^{-1})$ and $g_\ell:\star_{i=1}^{\ell+1}\Gamma(G_i^{-1})\rightarrow \star_{i=1}^{\ell}\Gamma(G_i^{-1})$ are defined by 
$$
f_{\ell}(x_1,x_2,x_3,\ldots, x_{\ell+1},x_{\ell+2})=(x_1,x_2,x_3,\ldots, x_{\ell+1})
$$
and
$$
g_{\ell}(y_1,y_2,y_3,\ldots, y_{\ell+1},y_{\ell+2})=(y_1,y_2,y_3,\ldots, y_{\ell+1})
$$
for any $(x_1,x_2,x_3,\ldots, ,x_{\ell+1},x_{\ell+2})\in \star_{i=1}^{\ell+1}\Gamma(F_i^{-1})$ and any $(y_1,y_2,y_3,\ldots, ,y_{\ell+1},y_{\ell+2})\in \star_{i=1}^{\ell+1}\Gamma(G_i^{-1})$. It follows from Theorem \ref{mappingTHM2}, that there is a continuous surjection from $\varprojlim\{\star_{i=1}^{\ell}\Gamma(F_i^{-1}),f_\ell\}_{\ell=1} ^\infty$ to $\varprojlim\{\star_{i=1}^{\ell}\Gamma(G_i^{-1}),g_\ell\}_{\ell=1} ^\infty$, if and only if there is a sequence $(H_\ell)$ of upper semicontinuous functions 
$$
H_\ell:\star_{i=1}^{\ell+1}\Gamma(F_i^{-1})\multimap \star_{i=1}^{\ell}\Gamma(G_i^{-1})
$$
such that 
\begin{enumerate}
\item[(a)]  for each positive integer $k$  and each  $(x_1,x_2,x_3,\ldots, x_{k+2},x_{k+3})\in \star_{i=1}^{k+2}\Gamma(F_i^{-1})$,
$$
f_{k}(H_{k+1}(x_1,x_2,x_3,\ldots,x_{k+2},x_{k+3}))\subseteq H_k(f_{k+1}(x_1,x_2,x_3,\ldots, x_{k+2},x_{k+3})),
$$
\item[(b)]  for each positive integer $j$, 
$$
\lim_{k\to\infty}\diam (f_{j,k}(H_{k}(P_{k+1}((x_1,x_2),(x_1,x_2,x_3),\ldots))))=0 
$$
for any $((x_1,x_2),(x_1,x_2,x_3),\ldots)\in \varprojlim\{\star_{i=1}^{\ell}\Gamma(F_i^{-1}),f_\ell\}$, where 
$$
P_{k+1}:\varprojlim\{\star_{i=1}^{\ell}\Gamma(F_i^{-1}),f_\ell\}\rightarrow \star_{i=1}^{k+1}\Gamma(F_i^{-1})
$$
is the standard $(k+1)$-st projection, 
\end{enumerate}
and this is equivalent to \ref{zwei2}.
\end{proof}

The rest of the paper is dedicated to the fixed point property of inverse limits with set-valued functions. 
The following theorem gives a characterization of inverse limits with set-valued bonding functions that do not have the fixed point property.
\begin{theorem}\label{tralala}
Let $\{X_\ell,F_\ell\}_{\ell=1}^{\infty}$ be an inverse sequence of compact metric spaces $X_\ell$ and upper semicontinuous set-valued functions $F_\ell:X_{\ell+1}\multimap X_\ell$.   

The following statements are equivalent.
\begin{enumerate}
\item The inverse limit $\varproj\{X_\ell,F_\ell\}_{\ell=1} ^\infty$ does not have the fixed point property.
\item  \label{zwei} There is a sequence $(H_\ell)$ of upper semicontinuous functions 
$$
H_\ell:\star_{i=1}^{\ell+1}\Gamma(F_i^{-1})\multimap \star_{i=1}^{\ell}\Gamma(F_i^{-1})
$$
such that 
\begin{enumerate}
\item[(a)] \label{jedom} for each positive integer $k$, for each  $(x_1,x_2,x_3,\ldots, x_{k+2},x_{k+3})\in \star_{i=1}^{k+2}\Gamma(F_i^{-1})$, and for each $(z_1,z_2,z_3,\ldots ,z_{k+1},z_{k+2})\in H_{k+1}(x_1,x_2,x_3,\ldots, x_{k+2},x_{k+3})$,
$$
(z_1,z_2,z_3,\ldots ,z_{k+1})\in H_{k}(x_1,x_2,x_3,\ldots, x_{k+2}).
$$
\item[(b)] \label{dvom} for each positive integer $j$ and for each $\mathbf x=(x_1,x_2,x_3,\ldots)\in \varproj\{X_\ell,F_\ell\}_{\ell=1} ^\infty$, 
$$
\lim_{k\to\infty}\diam (\mathcal{H}_{j,k}(\mathbf x))=0, 
$$
where for each positive integer $k>j$, $\mathcal{H}_{j,k}(\mathbf x)$ is the set of all points 
$$
(z_1,z_2,z_3,\ldots ,z_{j+1})\in \star_{i=1}^{j}\Gamma(F_i^{-1})
$$
such that there is a point $(y_1,y_2,y_3,\ldots,y_{k+1})\in H_k(x_1,x_2,x_3,\ldots, x_{k+2})$ such that $(y_1,y_2,y_3,\ldots,y_{j+1})=(z_1,z_2,z_3,\ldots,z_{j+1})$, and
\item[(c)] \label{trom} There is a positive integer $k$, such that 
$$
(x_1,x_2,x_3,\ldots, x_{k+1})\not \in H_k(x_1,x_2,x_3,\ldots, ,x_{k+1},x_{k+2})
$$
for each $(x_1,x_2,x_3,\ldots, ,x_{k+1}x_{k+2})\in \star_{i=1}^{k+1}\Gamma(F_i^{-1})$.
\end{enumerate}
\end{enumerate}
\end{theorem}
\begin{proof}
The inverse limit $\varproj\{X_\ell,F_\ell\}_{\ell=1} ^\infty$ does not have the fixed point property if and only if the inverse limit $\varprojlim\{\star_{i=1}^{\ell}\Gamma(F_i^{-1}),f_\ell\}_{\ell=1} ^\infty$ does not have the fixed point property, where each $f_\ell:\star_{i=1}^{\ell+1}\Gamma(F_i^{-1})\rightarrow \star_{i=1}^{\ell}\Gamma(F_i^{-1})$ is defined by 
$$
f_{\ell}(x_1,x_2,x_3,\ldots, ,x_{k+1},x_{k+2})=(x_1,x_2,x_3,\ldots, ,x_{k+1})
$$
for any $(x_1,x_2,x_3,\ldots, ,x_{k+1},x_{k+2})\in \star_{i=1}^{\ell+1}\Gamma(F_i^{-1})$. By Theorem \ref{main:3}, the inverse limit $\varprojlim\{\star_{i=1}^{\ell}\Gamma(F_i^{-1}),f_\ell\}_{\ell=1} ^\infty$ does not have the fixed point property if and only if there is a sequence $(H_\ell)$ of upper semicontinuous functions 
$$
H_\ell:\star_{i=1}^{\ell+1}\Gamma(F_i^{-1})\multimap \star_{i=1}^{\ell}\Gamma(F_i^{-1})
$$
such that 
\begin{enumerate}
\item[(a)]  for each positive integer $k$  and each  $(x_1,x_2,x_3,\ldots, x_{k+2},x_{k+3})\in \star_{i=1}^{k+2}\Gamma(F_i^{-1})$,
$$
f_{k}(H_{k+1}(x_1,x_2,x_3,\ldots,x_{k+2},x_{k+3}))\subseteq H_k(f_{k+1}(x_1,x_2,x_3,\ldots, x_{k+2},x_{k+3})),
$$
\item[(b)]  for each positive integer $j$, 
$$
\lim_{k\to\infty}\diam (f_{j,k}(H_{k}(P_{k+1}((x_1,x_2),(x_1,x_2,x_3),\ldots))))=0 
$$
for any $((x_1,x_2),(x_1,x_2,x_3),\ldots)\in \varprojlim\{\star_{i=1}^{\ell}\Gamma(F_i^{-1}),f_\ell\}$, where 
$$
P_{k+1}:\varprojlim\{\star_{i=1}^{\ell}\Gamma(F_i^{-1}),f_\ell\}\rightarrow \star_{i=1}^{k+1}\Gamma(F_i^{-1})
$$
is the standard $(k+1)$-st projection, and
\item[(c)]  there is a positive integer $k$, such that 
$$
f_k(x_1,x_2,x_3,\ldots, x_{k+1},x_{k+2})\not \in H_k(x_1,x_2,x_3,\ldots, ,x_{k+1},x_{k+2})
$$
for each $(x_1,x_2,x_3,\ldots, ,x_{k+1}x_{k+2})\in \star_{i=1}^{k+1}\Gamma(F_i^{-1})$,
\end{enumerate}
and this is equivalent to \ref{zwei}.
\end{proof}

We continue with the following theorems that only give sufficient conditions for an inverse limit not to have the fixed point property. However, these sufficient conditions are much more user friendly than the necessary and sufficient conditions that are obtained in Theorem \ref{tralala}. 
\begin{theorem}\label{Banic}
Let $\{X_\ell,F_\ell\}_{\ell=0}^{\infty}$ be an inverse sequence of compact metric spaces $X_\ell$ and upper semicontinuous set-valued functions $F_\ell:X_{\ell+1}\multimap X_\ell$, and let $\{X_\ell,g_\ell\}_{\ell=1}^{\infty}$ be an inverse sequence of compact metric spaces $X_\ell$ and continuous single-valued functions $g_\ell:X_{\ell+1}\rightarrow X_\ell$ such that for each positive integer $\ell$, 
$$
 g_\ell (F_{\ell+1}(x))\subseteq F_\ell(g_{\ell+1}(x))
$$
for any  $x \in X_{\ell+2}$.
 If there is an upper semicontinuous set-valued function $G_0:X_1\multimap X_0$ such that 
\begin{enumerate}
\item for any  $x \in X_{2}$, 
$$
G_0 (F_{1}(x)) \subseteq F_0(g_{1}(x)),
$$
\item   for any $x \in X_1$, 
$$
G_0(x) \not \subseteq F_0(x),
$$ 
\end{enumerate}
then $\varproj \{X_\ell,F_\ell\}_{\ell=1}^{\infty}$ does not have the fixed point property. [See Figure \ref{minioni}.]

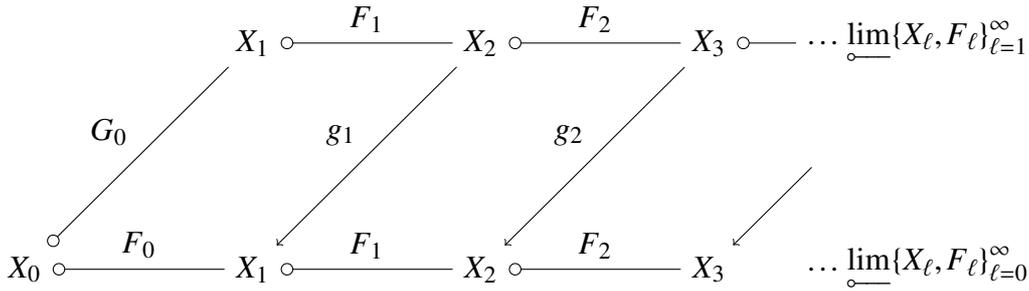
\begin{figure}[h]
\centering
\noindent \begin{tikzpicture}[node distance=1.5cm, auto]
  \node (X1) {};
  \node (X2) [right of=X1] {};
  \node (X3) [right of=X2] {$X_{1}$};
    \node (X4) [right of=X3] {};
  \node (X5) [right of=X4] {$X_{2}$};
   \node (X6) [right of=X5] {};
  \node (X7) [right of=X6] {$X_{3}$};
   \node (X8) [right of=X7] {$\ldots$};
    \node (X9) [right of=X8] {$\varproj\{X_\ell,F_\ell\}_{\ell=1}^{\infty}$};
  \draw[o-] (X3) to node {$F_{1}$} (X5);
   \draw[o-] (X5) to node {$F_{2}$} (X7);
    \draw[o-] (X7) to node {} (X8);
   \node (Z1) [below of=X1] {};
   \node (Z2) [right of=Z1] {};
   \node (Z3) [right of=Z2] {};
   \node (Z4) [right of=Z3] {};
   \node (Z5) [right of=Z4] {};
   \node (Z6) [right of=Z5] {};
   \node (Z7) [right of=Z6] {};
   \node (Z8) [right of=Z7] {};
   \node (Z9) [right of=Z8] {};
    \node (Y3) [below of=Z] {$X_{0}$};
          \node (Y4) [right of=Y3] {};
            \node (Y5) [right of=Y4] {$X_{1}$};
            \node (Y6) [right of=Y5] {};
            \node (Y7) [right of=Y6] {$X_{2}$};
            \node (Y8) [right of=Y7] {};
            \node (Y9) [right of=Y8] {$X_{3}$};
                 \node (Y10) [right of=Y9] {$\ldots$};
     \node (Y11) [right of=Y10] {$\varproj\{X_\ell,F_\ell\}_{\ell=0}^{\infty}$};
              \draw[o-] (Y3) to node {$F_{0}$} (Y5);
               \draw[o-] (Y5) to node {$F_{1}$} (Y7);
                \draw[o-] (Y7) to node {$F_{2}$} (Y9);
               \draw[o-] (Y3) to node {$G_0$} (X3);
                \draw[<-] (Y5) to node {$g_{1}$} (X5);
                 \draw[<-] (Y7) to node {$g_{2}$} (X7);
                 \draw[<-] (Y9) to node {} (Z8);
\end{tikzpicture}
\caption{A diagram from Theorem \ref{Banic}} \label{minioni}
\end{figure}
\end{theorem}

\begin{proof}
	Let $g:\varproj \{X_\ell,F_\ell\}_{\ell=1}^{\infty}\rightarrow \varproj \{X_\ell,F_\ell\}_{\ell=1}^{\infty}$ be defined by 
	$$
	g(x_1,x_2,x_3,\ldots)=(g_1(x_2), g_2(x_3), g_3(x_4),\ldots)
	$$
	for any $(x_1,x_2,x_3,\ldots) \in\varproj \{X_\ell,F_\ell\}_{\ell=1}^{\infty} $. We show that $g$ is a well-defined continuous function that does not have any fixed points. 
	
	Let $(x_1,x_2,x_3,\ldots) \in\varproj \{X_\ell,F_\ell\}_{\ell=1}^{\infty} $ be any point.  To show that $g(x_1,x_2,x_3,\ldots) \in\varproj \{X_\ell,F_\ell\}_{\ell=1}^{\infty} $, let $k$ be any positive integer.  We show that $g_k(x_{k+1})\in F_k(g_{k+1}(x_{k+2}))$. It follows from $x_{k+1}\in F_{k+1}(x_{k+2})$ and $g_k (F_{k+1}(x_{k+2}))\subseteq F_k(g_{k+1}(x_{k+2}))$ that
	$$
	g_k(x_{k+1})\in g_k(F_{k+1}(x_{k+2}))\subseteq F_k(g_{k+1}(x_{k+2})).
	$$
	We have just proved that $g$ is a well-defined function. Obviously, as such, it is a continuous function, since each $g_\ell$ is a continuous function. 
	
	Suppose next that there is a point $(x_1,x_2,x_3,\ldots) \in\varproj \{X_\ell,F_\ell\}_{\ell=1}^{\infty} $ such that 
	$$
	g(x_1,x_2,x_3,\ldots)=(g_1(x_2), g_2(x_3), g_3(x_4),\ldots)=(x_1,x_2,x_3,\ldots).
	$$
	Then $g_{\ell}(x_{\ell+1})=x_\ell$ for each positive integer $\ell$, and in the special case $\ell=1$ it means that $g_1(x_2)=x_1$ and we have
	$$
	F_0(x_1)=F_0(g_1(x_2))\supseteq G_0 (F_{1}(x_2))
	$$
	meaning that for each $x\in F_1(x_2)$,
	$$
	G_0(x) \subseteq F_0(x_1).
	$$
	holds.  Since $x_1\in F_1(x_2)$, it follows that
	$$
	G_0(x_1) \subseteq F_0(x_1),
	$$
	which is a contradiction with $G_0(x) \not \subseteq F_0(x)$ for any $x\in X_1$. Therefore,  for each $(x_1,x_2,x_3,\ldots) \in\varproj \{X_\ell,F_\ell\}_{\ell=1}^{\infty} $,
	$$
	g(x_1,x_2,x_3,\ldots)\neq (x_1,x_2,x_3,\ldots)
	$$
	holds.
\end{proof}
The following corollary is a well-known result that may also be obtained from the classic results from \cite{feuerbacher,mioduszewski}. It was used by many authors to construct an example of a tree-like continuum that does not have the fixed point property; for examples see \cite{logan,38,39}.
\begin{corollary}\label{Banic:cor}
Let $\{X_\ell,f_\ell\}_{\ell=0}^{\infty}$ be an inverse sequence of compact metric spaces $X_\ell$ and continuous  functions $f_\ell:X_{\ell+1}\rightarrow X_\ell$, and let $\{X_\ell,g_\ell\}_{\ell=0}^{\infty}$ be an inverse sequence of compact metric spaces $X_\ell$ and continuous single-valued functions $g_\ell:X_{\ell+1}\rightarrow X_\ell$ such that for each non-negative integer $\ell$, 
$$
 g_\ell (f_{\ell+1}(x)) = f_\ell(g_{\ell+1}(x))
$$
for any  $x \in X_{\ell+2}$.
 If $g_0(x) \neq f_0(x)$ for all $x\in X_1$,  then $\varprojlim \{X_\ell,f_\ell\}_{\ell=1}^{\infty}$ does not have the fixed point property. [See Figure \ref{minioni1}.]

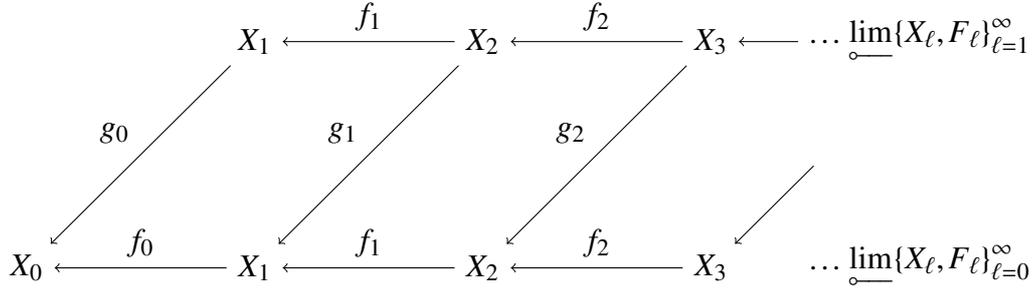
\begin{figure}[h]
\centering
\noindent \begin{tikzpicture}[node distance=1.5cm, auto]
  \node (X1) {};
  \node (X2) [right of=X1] {};
  \node (X3) [right of=X2] {$X_{1}$};
    \node (X4) [right of=X3] {};
  \node (X5) [right of=X4] {$X_{2}$};
   \node (X6) [right of=X5] {};
  \node (X7) [right of=X6] {$X_{3}$};
   \node (X8) [right of=X7] {$\ldots$};
    \node (X9) [right of=X8] {$\varproj\{X_\ell,F_\ell\}_{\ell=1}^{\infty}$};
  \draw[<-] (X3) to node {$f_{1}$} (X5);
   \draw[<-] (X5) to node {$f_{2}$} (X7);
    \draw[<-] (X7) to node {} (X8);
   \node (Z1) [below of=X1] {};
   \node (Z2) [right of=Z1] {};
   \node (Z3) [right of=Z2] {};
   \node (Z4) [right of=Z3] {};
   \node (Z5) [right of=Z4] {};
   \node (Z6) [right of=Z5] {};
   \node (Z7) [right of=Z6] {};
   \node (Z8) [right of=Z7] {};
   \node (Z9) [right of=Z8] {};
    \node (Y3) [below of=Z] {$X_{0}$};
          \node (Y4) [right of=Y3] {};
            \node (Y5) [right of=Y4] {$X_{1}$};
            \node (Y6) [right of=Y5] {};
            \node (Y7) [right of=Y6] {$X_{2}$};
            \node (Y8) [right of=Y7] {};
            \node (Y9) [right of=Y8] {$X_{3}$};
                 \node (Y10) [right of=Y9] {$\ldots$};
     \node (Y11) [right of=Y10] {$\varproj\{X_\ell,F_\ell\}_{\ell=0}^{\infty}$};
              \draw[<-] (Y3) to node {$f_{0}$} (Y5);
               \draw[<-] (Y5) to node {$f_{1}$} (Y7);
                \draw[<-] (Y7) to node {$f_{2}$} (Y9);
               \draw[<-] (Y3) to node {$g_0$} (X3);
                \draw[<-] (Y5) to node {$g_{1}$} (X5);
                 \draw[<-] (Y7) to node {$g_{2}$} (X7);
                 \draw[<-] (Y9) to node {} (Z8);
\end{tikzpicture}
\caption{A diagram from Corollary \ref{Banic:cor}} \label{minioni1}
\end{figure}
\end{corollary}

\begin{proof}
For each non-negative integer $\ell$, let $F_\ell:X_{\ell+1}\multimap X_\ell$ be defined by
$F_\ell(x)=\{f_\ell(x)\}$ for any $x\in X_{\ell+1}$, and let $G_0(x)=\{g_0(x)\}$ for any $x\in X_1$. Then 
\begin{enumerate}
\item for any  $x \in X_{2}$, 
$$
G_0 (F_{1}(x)) \subseteq F_0(g_{1}(x)),
$$
\item   for any $x \in X_1$, 
$$
G_0(x) \not \subseteq F_0(x).
$$ 
\end{enumerate}
Therefore, by Theorem \ref{Banic},  $\varproj \{X_\ell,F_\ell\}_{\ell=1}^{\infty}$ does not have the fixed point property.	It follows from  $\varproj \{X_\ell,F_\ell\}_{\ell=1}^{\infty}=\varprojlim \{X_\ell,f_\ell\}_{\ell=1}^{\infty}$ that the inverse limit $\varprojlim \{X_\ell,f_\ell\}_{\ell=1}^{\infty}$ does not have the fixed point property. 
\end{proof}

We conclude the paper with the following theorem.

\begin{theorem}
Let $\{X_\ell,F_\ell\}_{\ell=1}^{\infty}$ be an inverse sequence of compact metric spaces $X_\ell$ and upper semicontinuous set-valued functions $F_\ell:X_{\ell+1}\multimap X_\ell$, and let $(h_\ell)$ be a sequence of continuous functions  
$$
h_\ell:\star_{i=1}^{\ell+1}\Gamma(F_i^{-1})\rightarrow \star_{i=1}^{\ell}\Gamma(F_i^{-1})
$$ 
such that
\begin{enumerate}
\item $h_1(x_1,x_2,x_3)\neq (x_1,x_2)$ for each $(x_1,x_2,x_3)\in \star_{i=1}^{2}\Gamma(F_i^{-1})$, and
\item for each positive integer $k$  and each  $(x_1,x_2,x_3,\ldots, x_{k+2},x_{k+3})\in \star_{i=1}^{k+2}\Gamma(F_i^{-1})$,  
$$
(z_1,z_2,z_3,\ldots ,z_{k+1},z_{k+2})=h_{k+1}(x_1,x_2,x_3,\ldots, x_{k+2},x_{k+3}) 
$$
$$
\Longrightarrow (z_1,z_2,z_3,\ldots ,z_{k+1})=h_{k}(x_1,x_2,x_3,\ldots, x_{k+2}).
$$
\end{enumerate}
Then the inverse limit $\varproj\{X_\ell,F_\ell\}_{\ell=1} ^\infty $ does not have the fixed point property.
\end{theorem}
\begin{proof}
For each positive integer $k$, let 
\begin{enumerate}
\item $Y_{k-1}=\star_{i=1}^{k}\Gamma(F_i^{-1})$,
\item $f_{k-1}:Y_{k}\rightarrow Y_{k-1}$ be defined by 
$$
f_{k-1}(x_1,x_2,x_3,\ldots,x_{k+1},x_{k+2})=(x_1,x_2,x_3,\ldots,x_{k+1})
$$ 
for each $(x_1,x_2,x_3,\ldots,x_{k+1},x_{k+2})\in Y_{k}$, and 
\item $g_{k-1}=h_k$.
\end{enumerate}
 By Theorem \ref{vanizrek2}, the inverse limits 
$$
\varproj\{X_\ell,F_\ell\}_{\ell=1} ^\infty \textup{ ~~~ and ~~~ } \varprojlim\{Y_\ell, f_\ell\}_{n=0} ^\infty
$$ 
are homeomorphic.  It follows from Corollary  \ref{Banic:cor} that $\varprojlim \{Y_\ell,f_\ell\}_{\ell=1}^{\infty}$ does not have the fixed point property. Therefore, the inverse limit $\varproj\{X_\ell,F_\ell\}_{\ell=1} ^\infty$ does not have the fixed point property.
\end{proof}



\noindent I. Bani\v c\\
              (1) Faculty of Natural Sciences and Mathematics, University of Maribor, Koro\v{s}ka 160, SI-2000 Maribor,
   Slovenia; \\(2) Institute of Mathematics, Physics and Mechanics, Jadranska 19, SI-1000 Ljubljana, 
   Slovenia; \\(3) Andrej Maru\v si\v c Institute, University of Primorska, Muzejski trg 2, SI-6000 Koper,
   Slovenia\\
             {iztok.banic@um.si}           
     
				\-
				
		\noindent G.  Erceg\\
             Department of Mathematics, Faculty of Science, University of Split, Rudera Bo\v skovi\' ca 33, Split,  Croatia\\
{{gorerc@pmfst.hr}       }    

                 	\-
					
  \noindent J.  Kennedy\\
             Lamar University, 200 Lucas Building, P.O. Box 10047, Beaumont, TX 77710 USA\\
{{kennedy9905@gmail.com}       }    

\end{document}